\documentclass[a4paper]{amsart}
\usepackage{amscd, amssymb, color, booktabs, chemarrow}
\usepackage{myown}

\usepackage{harpoon}

\title{Finite Theta Correspondence of Almost Characters}
\author{Shu-Yen Pan}
\address{Department of Mathematics,
National Tsing Hua University, Hsinchu 300, Taiwan}
\email{sypan@math.nthu.edu.tw}


\date{\today}

\begin{document}

\begin{abstract}
The theory of almost characters which is closely related to character sheaves is proposed by Lusztig
to study the representation theory of finite reductive groups.
In this article we show that the decomposition of the Weil character for finite reductive dual pairs
$(\Sp_{2n},\rmO^\pm_{2n'})$ or $(\Sp_{2n},\SO_{2n'+1})$ with respect to the almost characters
is exactly the same as the decomposition with respect to the irreducible characters.
As a consequence, the finite theta correspondence on almost characters is established.
\end{abstract}

\maketitle
\tableofcontents

\section{Introduction}

\subsection{}
Let $\bfG$ be a classical group such as $\GL_n,\rmU_n,\Sp_{2n},\rmO^\epsilon_n$ defined
over a finite field $\bfF_q$ of odd characteristic,
and let $G$ denote the group of rational points.
The (complex) space $\calv(\bfG)$ of class functions on $G$ is an inner product space with
respect to the inner product
\[
\langle f,f'\rangle_\bfG=\frac{1}{|G|}\sum_{g\in G}f(g)f'(g^{-1})
\]
for $f,f'\in\calv(\bfG)$.
The set of irreducible characters
\[
\cale(\bfG)=\{\,\rho_x\mid x\in X(\bfG)\,\}
\]
of $G$ for some index set $X(\bfG)$ forms an orthonormal basis for $\calv(\bfG)$.
Another orthonormal basis
\[
\cala(\bfG)=\{\,R_x\mid x\in X(\bfG)\,\}
\]
whose elements are called \emph{almost characters} and are defined as linear combinations of
irreducible characters is proposed by Lusztig (\cf.~\cite{lusztig-book} p.347).
If $\bfG$ is a general linear group $\GL_n$ or a unitary group $\rmU_n$,
then $R_x=\rho_x$ for $x\in X(\bfG)$ and there is nothing really interesting.
However, if $\bfG$ is a symplectic group or an orthogonal group,
the story is very different.
For $\bfG=\Sp_{2n}$, $\rmO^\pm_{2n}$ or $\SO_{2n+1}$,
the precise definition of $\cala(\bfG)$ is given in \cite{waldspurger} \S 4.
Note that here we regard $\rmO^\pm$ as a ``group'' and define 
$\cale(\rmO^\pm_{2n})=\cale(\rmO^+_{2n})\cup\cale(\rmO^-_{2n})$,
$X(\rmO^\pm_{2n})=X(\rmO^+_{2n})\cup X(\rmO^-_{2n})$,
and then $\calv(\rmO^\pm_{2n})=\calv(\rmO^+_{2n})\oplus\calv(\rmO^-_{2n})$.

\subsection{}
If $(\bfG,\bfG')=(\Sp_{2n},\rmO^\epsilon_{2n'})$ or $(\Sp_{2n},\rmO_{2n'+1})$
where $\epsilon=+$ or $-$ is a reductive dual pair
of a symplectic group and an orthogonal group,
let $\omega^\psi_{\bfG,\bfG'}$ denote the Weil character of the dual pair $(\bfG,\bfG')$
with respect to a nontrivial additive character $\psi$ of $\bfF_q$.
We can also regard $(\Sp_{2n},\rmO^\pm_{2n'})$ as a ``dual pair'', 
and define
\begin{align*}
\omega^\psi_{\Sp_{2n},\rmO^\pm_{2n'}}
&=\omega^\psi_{\Sp_{2n},\rmO^+_{2n'}}+\omega^\psi_{\Sp_{2n},\rmO^-_{2n'}}\in\calv(\Sp_{2n})\otimes\calv(\rmO^\pm_{2n'}) \\
\omega^\psi_{\Sp_{2n},\SO_{2n'+1}}
&=\omega^\psi_{\Sp_{2n},\rmO_{2n'+1}}|_{\Sp_{2n}(q)\times\SO_{2n'+1}(q)}\in\calv(\Sp_{2n})\otimes\calv(\SO_{2n'+1}).
\end{align*}
The decomposition of $\omega^\psi_{\Sp_{2n},\rmO^\pm_{2n'}}$ (resp.~$\omega^\psi_{\Sp_{2n},\SO_{2n'+1}}$)
with respect to the orthonormal basis $\cale(\Sp_{2n}\times\rmO^\pm_{2n'})$
(resp.~$\cale(\Sp_{2n}\times\SO_{2n'+1})$) consisting of irreducible characters is known in
\cite{pan-finite-unipotent}, \cite{pan-Lusztig-correspondence}.
The following main result of this article is to show that they have the same decomposition with respect to
the orthonormal basis consisting of almost characters $\cala(\Sp_{2n}\times\rmO^\pm_{2n'})$
(resp.~$\cala(\Sp_{2n}\times\SO_{2n'+1})$).

\begin{thm}\label{0101}
Let $(\bfG,\bfG')=(\Sp_{2n},\rmO^\pm_{2n'})$ or $(\Sp_{2n},\SO_{2n'+1})$,
and write
\begin{align*}
\omega^\psi_{\bfG,\bfG'}
&=\sum_{x\in X(\bfG),\ x'\in X(\bfG')} m_{x,x'}\rho_x\otimes\rho_{x'} \\
&=\sum_{x\in X(\bfG),\ x'\in X(\bfG')} \widetilde m_{x,x'}R_x\otimes R_{x'}
\end{align*}
for $m_{x,x'},\widetilde m_{x,x'}\in\bbC$.
Then $m_{x,x'}=\widetilde m_{x,x'}$ for any $(x,x')\in X(\bfG)\times X(\bfG')$.
\end{thm}

It might be a little surprised that the above result does not hold for the dual pair
$(\Sp_{2n},\rmO^\epsilon_{2n'})$ or $(\Sp_{2n},\rmO_{2n'+1})$.
An example is given in Subsection~\ref{0326}.
From this point of view, the dual pairs $(\Sp_{2n},\rmO^\pm_{2n'})$ and $(\Sp_{2n},\SO_{2n'+1})$
might be more natural than the usual $(\Sp_{2n},\rmO^+_{2n'})$, $(\Sp_{2n},\rmO^-_{2n'})$ or
$(\Sp_{2n},\rmO_{2n'+1})$.

It is well known that the correspondence for dual pair $(\Sp_{2n},\rmO^\epsilon_{2n'})$ or
$(\Sp_{2n},\rmO_{2n'+1})$ is of multiplicity one.
Therefore, the correspondence for the pair $(\Sp_{2n},\rmO^\pm_{2n'})$ is also of multiplicity one,
i.e., the coefficient $m_{x,x'}$ in Theorem~\ref{0101} is either $0$ or $1$.
Although the correspondence for $(\Sp_{2n},\SO^\epsilon_{2n'})$ is not of multiplicity one
(\cf.~Example~\ref{0423}), the correspondence for the pair $(\Sp_{2n},\SO_{2n'+1})$ is still of multiplicity one.
This result should be known to the experts, still we provide a proof (\cf.~Proposition~\ref{0424}) for completeness.

A direct consequence of the above theorem is that one can define a correspondence on almost
characters for dual pair $(\Sp_{2n},\rmO^\pm_{2n'})$ or $(\Sp_{2n},\SO_{2n'+1})$
which is exactly the same as the correspondence on irreducible characters.
As an application, we can define the $\Theta$-rank of an almost character $R_x$
as we did for an irreducible character $\rho_x$.
Moreover, the above theorem tells us that $\Theta\text{\rm -rk}(R_x)=\Theta\text{\rm -rk}(\rho_x)$
for any $x\in X(\bfG)$ for $\bfG=\Sp_{2n},\rmO^\pm_{2n}$ or $\SO_{2n+1}$.

\subsection{}
Now we briefly explain the idea of our proof.
It is known that the set of irreducible characters
\[
\cale(\bfG)=\bigcup_{(s)\subset (G^*)^0\,\text{semisimple}}\cale(\bfG)_s
\]
is partition into \emph{Lusztig series} indexed by conjugacy classes of semisimple
elements in the connected component of the dual group.
Elements in $\cale(\bfG)_1$ are called \emph{unipotent} and parametrized by ``\emph{symbols}'' $\cals_\bfG$.
The set of symbols are partitioned into families $\cals_\bfG=\bigcup_Z\cals_Z$ indexed by the set of
``\emph{special symbols}''.

For each family of symbols $\cals_Z$, the transition matrix between irreducible characters
and almost characters is known explicitly in \cite{lusztig-book} or \cite{waldspurger}.
Then we can adapt the technique from \cite{pan-uniform} to prove the unipotent part of the
equality of Theorem~\ref{0101}:
\begin{equation}\label{0104}
\omega_{\Sp_{2n},\rmO^\pm_{2n'},1}
=\sum_{(\Lambda,\Lambda')\in\calb_{\Sp_{2n},\rmO^\pm_{2n'}}}\rho_\Lambda\otimes\rho_{\Lambda'}
=\sum_{(\Lambda,\Lambda')\in\calb_{\Sp_{2n},\rmO^\pm_{2n'}}}R_\Lambda\otimes R_{\Lambda'}
\end{equation}
(\cf.~Corollary~\ref{0306}).

For a general Lusztig series $\cale(\bfG)_s$, we can define groups
$\bfG^{(0)}(s),\bfG^{(-)}(s),\bfG^{(+)}(s)$ and a \emph{modified Lusztig correspondence}
\[
\Xi_s\colon\cale(\bfG)_s\rightarrow \cale(\bfG^{(0)}(s)\times\bfG^{(-)}(s)\times\bfG^{(+)}(s))_1
\]
which is a bijection and can be extended linearly to an isometry of inner product spaces
\[
\Xi_s\colon\calv(\bfG)_s\rightarrow \calv(\bfG^{(0)}(s))_1\otimes\calv(\bfG^{(-)}(s))_1\otimes\calv(\bfG^{(+)}(s))_1.
\]
Then we can combine several Lusztig series together so that $\bfG^{(-)}(s),\bfG^{(+)}(s)$
are enlarged to $\bfG^{(-)}[s],\bfG^{(+)}[s]$ and each of which is of the form
$\Sp_{\even},\rmO^\pm_{\even}$ or $\SO_{\odd}$.
Then we have an isometry of inner product spaces
\[
\Xi_{[s]}\colon\calv(\bfG)_{[s]}
\rightarrow \calv(\bfG^{(0)}[s])_1\otimes\calv(\bfG^{(-)}[s])_1\otimes\calv(\bfG^{(+)}[s])_1.
\]
Now the almost character structure of $\calv(\bfG^{(0)}[s])_1\otimes\calv(\bfG^{(-)}[s])_1\otimes\calv(\bfG^{(+)}[s])_1$
induces the almost character structure of $\calv(\bfG)_{[s]}$ via the isometry $\Xi_{[s]}$.
If particular, if $\Xi_{[s]}(\rho_x)=\rho_{x^{(0)}}\otimes\rho_{\Lambda_x^{(-)}}\otimes\rho_{\Lambda_x^{(+)}}$,
then $\Xi_{[s]}(R_x)=R_{x^{(0)}}\otimes R_{\Lambda_x^{(-)}}\otimes R_{\Lambda_x^{(+)}}$.
Via the isometry $\Xi_{[s]}$, the identity of unipotent class functions in (\ref{0104})
can be extend to an identity of all class functions for $(\bfG,\bfG')=(\Sp_{2n},\rmO^\pm_{2n'})$
or $(\Sp_{2n},\SO_{2n'+1})$.
Then Theorem~\ref{0101} is obtained.

\subsection{}
The contents of this article are as follows.
In Section 2 we recall the basic notation of symbols by Lusztig.
Then we provide the definition of unipotent almost characters $\cala(\bfG)_1$ for $\bfG=\Sp_{2n}$ or
$\rmO^\pm_{2n}$ from \cite{waldspurger}.
In Section 3, we recall the result on the theta correspondence on unipotent characters from
\cite{pan-finite-unipotent} and then prove (\ref{0104}) with the technique adapted from \cite{pan-uniform}.
In Section 4, we prove our main results using modified Lusztig correspondences from
\cite{pan-Lusztig-correspondence}.
In the final section, we provide an application of our main result.

\section{Symbols and Unipotent Characters}

\subsection{Bi-partitions and symbols}
An ordered pair $\sqbinom{\lambda}{\mu}$ of two partitions $\lambda,\mu$ is called a \emph{bi-partition}.
For a non-negative integer $n$,
let $\calp_2(n)$ denote the set of bi-partitions of $n$, i.e.,
\[
\calp_2(n)=\left\{\,\textstyle\sqbinom{\lambda}{\mu}\mid |\lambda|+|\mu|=n\,\right\}.
\]

A \emph{symbol} $\Lambda$ is an ordered pair
\[
\Lambda=\binom{A}{B}=\binom{a_1,a_2,\ldots,a_{m_1}}{b_1,b_2,\ldots,b_{m_2}}
\]
of two finite subsets $A,B$ (possibly empty) of non-negative integers whose
elements are always written in decreasing order, i.e., $a_1>a_2>\cdots>a_{m_1}$ and
$b_1>b_2>\cdots>b_{m_2}$.
For a symbol $\Lambda=\binom{A}{B}$, let $\Lambda^*=A$ and $\Lambda_*=B$
denote the first row and the second row of $\Lambda$ respectively.
Moreover, we define its \emph{transpose} $\Lambda^\rmt=\binom{B}{A}$,
its number of entries, \emph{size}, \emph{rank}, and \emph{defect} by:
\begin{align*}
|\Lambda| &=|A|+|B|, \\
{\rm size}(\Lambda) &= (|A|,|B|), \\
{\rm rk}(\Lambda)
&=\sum_{i=1}^{m_1}a_i+\sum_{j=1}^{m_2}b_j-\left\lfloor\biggl(\frac{|A|+|B|-1}{2}\biggr)^2\right\rfloor, \\
{\rm def}(\Lambda) &=|A|-|B|.
\end{align*}

We define an equivalence relation on symbols generated by
\[
\binom{a_1,a_2,\ldots,a_{m_1}}{b_1,b_2,\ldots,b_{m_2}}\sim
\binom{a_1+1,a_2+1,\ldots,a_{m_1}+1,0}{b_1+1,b_2+1,\ldots,b_{m_2}+1,0},
\]
and let $\cals$ (resp.~$\cals_{n,\delta}$) denote the set of equivalence classes of symbols
(resp.~equivalences of symbols of rank $n$ and defect $\delta$).
It is clear that each equivalence class of symbols contains a unique symbol $\binom{A}{B}$ such that
$0\not\in A\cap B$.
Such a symbol is called \emph{reduced}.
Unless otherwise specified, a symbol is always assumed to be reduced.
It is not difficult to check that the mapping
\begin{equation}\label{0209}
\Upsilon\colon\binom{a_1,a_2,\ldots,a_{m_1}}{b_1,b_2,\ldots,b_{m_2}}\mapsto
\sqbinom{a_1-(m_1-1),a_2-(m_1-2),\cdots,a_{m_1-1}-1,a_{m_1}}{b_1-(m_2-1),b_2-(m_2-2),\cdots,b_{m_2-1}-1,b_{m_2}}
\end{equation}
gives a bijection $\cals_{n,1}\rightarrow\calp_2(n)$ and a bijection $\cals_{n,0}\rightarrow\calp_2(n)$.
Note that $\Upsilon(\Lambda_1)=\Upsilon(\Lambda_2)$ if $\Lambda_1\sim\Lambda_2$.

\subsection{Unipotent almost characters of a symplectic group}
Define
\[
\cals_{\Sp_{2n}}
=\bigcup_{\delta\equiv 1\pmod 4}\cals_{n,\delta}
=\{\,\Lambda\in\cals\mid{\rm rk}(\Lambda)=n,\ {\rm def}(\Lambda)\equiv 1\pmod 4\,\}
\]
It is known that there exists a parametrization of the unipotent characters $\cale(\Sp_{2n})_1$
of $\Sp_{2n}(q)$ by $\cals_{\Sp_{2n}}$ (\cf.~\cite{lg}).
For $\Lambda\in\cals_{\Sp_{2n}}$, let $\rho_\Lambda\in\cale(\Sp_{2n})_1$ denote the irreducible
character parametrized by $\Lambda$.
Then let $\calv(\Sp_{2n})_1$ denote the subspace spanned by $\{\,\rho_\Lambda\mid\Lambda\in\cals_{\Sp_{2n}}\,\}$.

A symbol
$Z=\binom{a_1,a_2,\ldots,a_{m+1}}{b_1,b_2,\ldots,b_m}$ of rank $n$ and defect $1$
is called \emph{special} if
$a_1\geq b_1\geq a_2\geq b_2\geq\cdots\geq a_m\geq b_m\geq a_{m+1}$.
For a special symbol $Z$, let $Z_\rmI$ denote the subsymbol of ``\emph{singles}'' in $Z$, i.e.,
$Z_\rmI=Z\smallsetminus\binom{Z^*\cap Z_*}{Z^*\cap Z_*}$.
Then $|Z_\rmI|=2d+1$ for some non-negative integer $d$,
called the \emph{degree} of $Z$, and denoted by $\deg(Z)=d$.
Let $\cals_Z$ denote the set of (reduced) symbols with exactly the same entries of $Z$, i.e.,
\[
\cals_Z=\{\,\Lambda\in\cals_{\Sp_{2n}}\mid \Lambda^*\cup\Lambda_*=Z^*\cup  Z_*,\ \Lambda^*\cap\Lambda_*=Z^*\cap Z_*\,\}.
\]
Then we have
\begin{align*}
\cals_{\Sp_{2n}} &=\bigcup_{Z\ \text{special},\ {\rm rk}(Z)=n,\ {\rm def}(Z)=1}\cals_Z, \\
\calv(\Sp_{2n})_1 &=\bigoplus_{Z\ \text{special},\ {\rm rk}(Z)=n,\ {\rm def}(Z)=1}\calv_Z
\end{align*}
where $\calv_Z=\calv(\Sp_{2n})_Z$ is the subspace of class functions spanned by
$\{\,\rho_\Lambda\mid\Lambda\in\cals_Z\,\}$.

For a subsymbol $M\subset Z_\rmI$,
we define
\begin{equation}\label{0210}
\Lambda_M=(Z\smallsetminus M)\cup M^\rmt,
\end{equation}
i.e., $\Lambda_M$ is obtained from $Z$ by switching the two rows of $M$.
Clearly, ${\rm rk}(\Lambda_M)={\rm rk}(Z)$.
Then it is easy to see that
\[
\cals_Z=\{\,\Lambda_M\mid M\subset Z_\rmI,\ |M|\text{ even}\,\},
\]
and $|\cals_Z|=2^{2\deg(Z)}$.
For $\Lambda_{M_1},\Lambda_{M_2}\in\cals_Z$, we define an addition
\begin{equation}\label{0211}
\Lambda_{M_1}+\Lambda_{M_2}=\Lambda_N\quad\text{where }N=(M_1\cup M_2)\smallsetminus(M_1\cap M_2).
\end{equation}
This gives $\cals_Z$ a vector space structure over $\bfF_2$ with identity element $\Lambda_{\binom{-}{-}}=Z$.
Moreover, we have a symmetric pairing $\langle,\rangle\colon\cals_Z\times\cals_Z\rightarrow\bfF_2$ given by
\begin{equation}\label{0203}
\langle\Lambda_{M_1},\Lambda_{M_2}\rangle=|M_1\cap M_2|\pmod 2.
\end{equation}

According to Lusztig, the \emph{unipotent almost character} $R_\Lambda$ of $\Sp_{2n}(q)$
associated to $\Lambda\in\cals_Z\subset\cals_{\Sp_{2n}}$ is defined by
\begin{equation}\label{0205}
R_\Lambda=\frac{1}{2^{\deg(Z)}}\sum_{\Lambda_1\in\cals_Z}(-1)^{\langle\Lambda,\Lambda_1\rangle}\rho_{\Lambda_1}
\in\calv(\Sp_{2n})_Z\subset\calv(\Sp_{2n})_1.
\end{equation}
This means that
\[
\langle R_{\Lambda_1},\rho_{\Lambda_2}\rangle_{\Sp_{2n}}
=\begin{cases}
\frac{1}{2^{\deg(Z)}}(-1)^{\langle\Lambda_1,\Lambda_2\rangle}, 
& \text{if $\Lambda_1,\Lambda_2\in\cals_Z$ for some special $Z$};\\
0, & \text{otherwise}.
\end{cases}
\]
In particular, if $\deg(Z)=0$, then $\cals_Z=\{Z\}$ and $R_Z=\rho_Z$.
The set of unipotent almost characters of $\Sp_{2n}(q)$ is denoted by
\[
\cala(\Sp_{2n})_1=\{\, R_\Lambda\mid\Lambda\in\cals_{\Sp_{2n}}\,\}.
\]

\begin{lem}\label{0207}
Let $Z$ be a special symbol of defect $1$.
The set $\cala_Z:=\{\,R_\Lambda\mid\Lambda\in\cals_Z\,\}$ forms an orthonormal basis of $\calv_Z$.
\end{lem}
\begin{proof}
Let $M_1,M_2\subset Z_\rmI$ such that both $|M_1|$ and $|M_2|$ are even.
Then by (\ref{0203}), we have
\begin{align*}
\langle R_{\Lambda_{M_1}},R_{\Lambda_{M_2}}\rangle
&= \frac{1}{2^{2\deg(Z)}}\sum_{M\subset Z_\rmI,\ |M|\ \text{even}}\;\sum_{M'\subset Z_\rmI,\ |M'|\ \text{even}}
(-1)^{|M\cap M_2|+|M'\cap M_2|}\langle\rho_{\Lambda_M},\rho_{\Lambda_{M'}}\rangle \\
&= \frac{1}{2^{2\deg(Z)}}\sum_{M\subset Z_\rmI,\ |M|\ \text{even}}(-1)^{|M\cap M_1|+|M\cap M_2|}
\end{align*}
If $M_1=M_2$, then clearly
\[
\sum_{M\subset Z_\rmI,\ |M|\ \text{even}}(-1)^{|M\cap M_1|+|M\cap M_2|}=|\cals_Z|=2^{2\deg(Z)}.
\]
Suppose that $M_1\neq M_2$.
Without loss of generality, may assume that there is an entry $a_1\in M_2\smallsetminus M_1$
and another entry $a_2\in Z_\rmI\smallsetminus M_2$ (because $M_2\subsetneq Z_\rmI$).
It is not difficult to see that
\[
\{\, M\mid M\subset Z_\rmI,\ |M|\ \text{even}\,\}
=\{\, M, M^*\mid M\subset Z_\rmI,\ |M|\ \text{even},\ a_1\not\in M\}
\]
where
\[
M^*=\begin{cases}
M\cup\{a_1,a_2\}, & \text{if $a_2\not\in M$};\\
(M\cup\{a_1\})\smallsetminus\{a_2\}, & \text{if $a_2\in M$}.
\end{cases}
\]
Now
\[
|M^*\cap M_1|+|M^*\cap M_2|=|M\cap M_1|+|M\cap M_2|+1
\]
for any $M\subset Z_\rmI$ such that $|M|$ even and $a_1\not\in M$.
Hence, if $M_1\neq M_2$, we have
\[
\sum_{M\subset Z_\rmI,\ |M|\ \text{even}}(-1)^{|M\cap M_1|+|M\cap M_2|}=0.
\]
Therefore the lemma is proved.
\end{proof}

\begin{exam}
The group $\bfG=\Sp_4$ has $6$ unipotent characters which are divided into three families
corresponding the special symbols $\binom{2}{-},\binom{2,0}{1},\binom{2,1,0}{2,1}$ respectively.
\begin{enumerate}
\item For $Z=\binom{2}{-}$, the family $\cale(\Sp_4)_Z$ has a unique character $\rho_{\binom{2}{-}}={\bf1}_{\Sp_4}$,
the trivial character of $\Sp_4(q)$, and then $R_{\binom{2}{-}}=\rho_{\binom{2}{-}}$.

\item For $Z=\binom{2,0}{1}$, the family $\cale(\Sp_4)_Z$ has $4$ characters
$\rho_{\binom{2,0}{1}},\rho_{\binom{2,1}{0}},\rho_{\binom{1,0}{2}},\rho_{\binom{-}{2,1,0}}$.
Then the table of $(-1)^{\langle\Lambda_1,\Lambda_2\rangle}$ for $\Lambda_1,\Lambda_2\in\cals_Z$ is:
\[
\begin{tabular}{ccc|cccc}
\toprule
& & $\Lambda_{M_2}$ & $\binom{2,0}{1}$ & $\binom{2,1}{0}$ & $\binom{1,0}{2}$ & $\binom{-}{2,1,0}$ \\
& & $M_2$ & $\binom{-}{-}$ & $\binom{0}{1}$ & $\binom{2}{1}$ & $\binom{2,0}{-}$ \\
$\Lambda_{M_1}$ & $M_1$ & & & & & \\
\midrule
$\binom{2,0}{1}$ & $\binom{-}{-}$ & & $\phantom{-}1$ & $\phantom{-}1$ & $\phantom{-}1$ & $\phantom{-}1$ \\
$\binom{2,1}{0}$ & $\binom{0}{1}$ & & $\phantom{-}1$ & $\phantom{-}1$ & $-1$ & $-1$ \\
$\binom{1,0}{2}$ & $\binom{2}{1}$ & & $\phantom{-}1$ & $-1$ & $\phantom{-}1$ & $-1$ \\
$\binom{-}{2,1,0}$ & $\binom{2,0}{-}$ & & $\phantom{-}1$ & $-1$ & $-1$ & $\phantom{-}1$ \\
\bottomrule
\end{tabular}
\]
Therefore,
\begin{align*}
R_{\binom{2,0}{1}}
&=\frac{1}{2}\bigl[\rho_{\binom{2,0}{1}}+\rho_{\binom{2,1}{0}}+\rho_{\binom{1,0}{2}}+\rho_{\binom{-}{2,1,0}}\bigr], \\
R_{\binom{2,1}{0}}
&=\frac{1}{2}\bigl[\rho_{\binom{2,0}{1}}+\rho_{\binom{2,1}{0}}-\rho_{\binom{1,0}{2}}-\rho_{\binom{-}{2,1,0}}\bigr], \\
R_{\binom{1,0}{2}}
&=\frac{1}{2}\bigl[\rho_{\binom{2,0}{1}}-\rho_{\binom{2,1}{0}}+\rho_{\binom{1,0}{2}}-\rho_{\binom{-}{2,1,0}}\bigr], \\
R_{\binom{-}{2,1,0}}
&=\frac{1}{2}\bigl[\rho_{\binom{2,0}{1}}-\rho_{\binom{2,1}{0}}-\rho_{\binom{1,0}{2}}+\rho_{\binom{-}{2,1,0}}\bigr].
\end{align*}

\item For $Z=\binom{2,1,0}{2,1}$, the family $\cale(\Sp_4)_Z$ has a unique character $\rho_{\binom{2}{-}}$,
the Steinberg character of $\Sp_4(q)$, and then $R_{\binom{2,1,0}{2,1}}=\rho_{\binom{2,1,0}{2,1}}$.
\end{enumerate}
\end{exam}

Let $\cals_{Z,1}$ be the set of symbols of defect $1$ in $\cals_Z$.
If $\Lambda\in\cals_{Z,1}$, then $\Upsilon(\Lambda)\in\calp_2(n)$ (\cf.~(\ref{0209})).
It is known that the set $\cale(W_n)$ of irreducible characters of the Weyl group $W_n$ of $\Sp_{2n}$
is parametrized by $\calp_2(n)$.
For $\Lambda\in\cals_{Z,1}\subset\cals_{n,1}$,
let $\chi_\Lambda\in\cale(W_n)$ be the irreducible character parametrized by $\Upsilon(\Lambda)$.
Then it is known by Lusztig (\cf.~\cite{pan-uniform} section 3) that 
\[
R_\Lambda=\frac{1}{|W_n|}\sum_{w\in W_n}\chi_\Lambda(w)R_{\bfT_w,1}
\]
where $R_{\bfT_w,1}=R_{\bfT_w,1}^{\Sp_{2n}}$ denotes the Deligne-Lusztig virtual character of $\Sp_{2n}(q)$
with respect to the trivial character of the rational maximal torus $T_w$.

\subsection{Unipotent almost characters of an even orthogonal group}
Now we consider the unipotent characters of an even orthogonal group
$\bfG=\rmO^+_{2n}$ or $\rmO^-_{2n}$.
Define
\begin{align*}
\cals_{\rmO^+_{2n}} &=\{\,\Lambda\mid{\rm rk}(\Lambda)=n,\ {\rm def}(\Lambda)\equiv 0\pmod 4\,\}, \\
\cals_{\rmO^-_{2n}} &=\{\,\Lambda\mid{\rm rk}(\Lambda)=n,\ {\rm def}(\Lambda)\equiv 2\pmod 4\,\}, \\
\cals_{\rmO^\pm_{2n}} := \cals_{\rmO^+_{2n}}\cup\cals_{\rmO^-_{2n}}
&=\{\,\Lambda\mid{\rm rk}(\Lambda)=n,\ {\rm def}(\Lambda)\equiv 0\pmod 2\,\}.
\end{align*}
Recall that $\cale(\rmO^\pm_{2n})_1:=\cale(\rmO^+_{2n})_1\cup\cale(\rmO^-_{2n})_1$.
Then there is a parametrization of $\cale(\rmO^\pm_{2n})_1$ by $\cals_{\rmO^\pm_{2n}}$.
For $\Lambda\in\cals_{\rmO^\pm_{2n}}$, let $\rho_\Lambda\in\cale(\rmO^\pm_{2n})_1$ denote the irreducible
character parametrized by $\Lambda$.
Then let $\calv(\rmO^\pm_{2n})_1$ denote the subspace spanned by
$\{\,\rho_\Lambda\mid\Lambda\in\cals_{\rmO^\pm_{2n}}\,\}$.

A symbol $Z=\binom{a_1,a_2,\ldots,a_m}{b_1,b_2,\ldots,b_m}$ of rank $n$ and defect $0$ is called 
\emph{special} if $a_1\geq b_1\geq a_2\geq b_2\geq\cdots\geq b_{m-1}\geq a_m\geq b_m$.
Again, let $Z_\rmI$ denote the subsymbol of ``\emph{singles}''.
Then $|Z_\rmI|=2d$ for some non-negative integer $d$,
called the \emph{degree} of $Z$.
Let $\cals_Z\subset\cals_{\rmO^\pm_{2n}}$ denote the set of symbols with exactly the same entries of $Z$.
Then
\begin{align*}
\cals_{\rmO^\pm_{2n}} &=\bigcup_{Z\ \text{special},\ {\rm rk}(Z)=n,\ {\rm def}(Z)=0}\cals_Z, \\
\calv(\rmO^\pm_{2n})_1 &=\bigoplus_{Z\ \text{special},\ {\rm rk}(Z)=n,\ {\rm def}(Z)=0}\calv_Z
\end{align*}
where $\calv_Z=\calv(\rmO^\pm_{2n})_Z$ denotes the space of class functions spanned by
$\{\,\rho_\Lambda\mid\Lambda\in\cals_Z\}$.

For a subsymbol $M\subset Z_\rmI$, as in (\ref{0210}), we define
\begin{equation}\label{0202}
\Lambda_M=(Z\smallsetminus M)\cup M^\rmt.
\end{equation}
It is clear that ${\rm rk}(\Lambda_M)={\rm rk}(Z)$,
${\rm def}(\Lambda_M)\equiv 0\pmod 4$ if $|M|$ is even;
and ${\rm def}(\Lambda_M)\equiv 2\pmod 4$ if $|M|$ is odd.
Moreover, we have
\begin{align*}
\cals_Z &=\{\,\Lambda_M\mid M\subset Z_\rmI\,\}, \\
\cals^+_Z:=\cals_Z\cap\cals_{\rmO^+_{2n}} &=\{\,\Lambda_M\mid M\subset Z_\rmI,\ |M|\text{ even}\,\}, \\
\cals^-_Z:=\cals_Z\cap\cals_{\rmO^-_{2n}} &=\{\,\Lambda_M\mid M\subset Z_\rmI,\ |M|\text{ odd}\,\},
\end{align*}
and $|\cals_Z|=2^{2\deg(Z)}$, $|\cals_Z^+|=|\cals_Z^-|=2^{2\deg(Z)-1}$.

For $\Lambda_{M_1},\Lambda_{M_2}\in\cals_Z$, we define $\Lambda_{M_1}+\Lambda_{M_2}$ 
as in (\ref{0211}).
This gives $\cals_Z$ a vector space structure over $\bfF_2$.
Note that if $\Lambda_1\in\cals_Z^{\epsilon_1}$ and $\Lambda_2\in\cals^{\epsilon_2}_Z$,
then $\Lambda_1+\Lambda_2\in\cals^{\epsilon_1\epsilon_2}_Z$ for $\epsilon_i=+$ or $-$.
Define a symmetric pairing $\langle,\rangle\colon\cals_Z\times\cals_Z\rightarrow\bfF_2$ by
\begin{equation}\label{0204}
\langle\Lambda_{M_1},\Lambda_{M_2}\rangle=|M_1||M_2|+|M_1\cap M_2|\pmod 2.
\end{equation}
Note that the definitions in (\ref{0203}) and (\ref{0204}) are consistent,
because $|M_i|$ is always assumed to be even in (\ref{0203}).
The following lemma says that the definition above is the same as that given in \cite{waldspurger} \S 2.1:

\begin{lem}
Let $Z$ be a special symbol of defect $0$.
For $\Lambda_1,\Lambda_2\in\cals_Z$, we have
\[
\langle\Lambda_1,\Lambda_2\rangle=(|(\Lambda_1)^*|+|Z^*|)(|(\Lambda_2)^*|+|Z^*|)
+|(\Lambda_1)^*\cap(\Lambda_2)^*\cap Z_*|+|(\Lambda_1)_*\cap(\Lambda_2)_*\cap Z^*|\pmod 2.
\]
\end{lem}
\begin{proof}
Suppose that $\Lambda_1=\Lambda_{M_1}$ and $\Lambda_2=\Lambda_{M_2}$ for some $M_1,M_2\subset Z_\rmI$.
For $M\subset Z_\rmI$, if $|M|$ is even, then $|(\Lambda_M)^*|\equiv |Z^*|\pmod 2$;
if $|M|$ is odd, then $|(\Lambda_M)^*|\not\equiv |Z^*|\pmod 2$.
Hence we have
\[
(|(\Lambda_1)^*|+|Z^*|)(|(\Lambda_2)^*|+|Z^*|)\equiv |M_1||M_2|\pmod 2.
\]
Moreover, it is easy to check that
\begin{align*}
(\Lambda_1)^*\cap(\Lambda_2)^*\cap Z_* &= (M_1)_*\cap(M_2)_*\cup(Z^*\cap Z_*) \\
(\Lambda_1)_*\cap(\Lambda_2)_*\cap Z^* &= (M_1)^*\cap(M_2)^*\cup(Z^*\cap Z_*) \\
\emptyset &= (M_1)_*\cap(M_2)^* =(M_1)^*\cap(M_2)_*.
\end{align*}
Therefore, we have
\[
|(\Lambda_1)^*\cap(\Lambda_2)^*\cap Z_*|+|(\Lambda_1)_*\cap(\Lambda_2)_*\cap Z^*|\equiv|M_1\cap M_2|\pmod 2.
\]
Hence, the lemma is obtained from (\ref{0204}).
\end{proof}

Following \cite{waldspurger} section 4,
the \emph{unipotent almost character} $R_\Lambda$ of $\rmO^\pm_{2n}(q)$
associated to a symbol $\Lambda\in\cals_Z$ is defined by
\begin{equation}\label{0206}
R_\Lambda
=\frac{1}{2^{\deg(Z)}}\sum_{\Lambda_1\in\cals_Z}(-1)^{\langle\Lambda,\Lambda_1\rangle}
\rho_{\Lambda_1}\in\calv(\rmO^\pm_{2n})_Z\subset\calv(\rmO^\pm_{2n})_1.
\end{equation}
Similarly, we have
\[
\langle R_{\Lambda_1},\rho_{\Lambda_2}\rangle_{\rmO^\pm_{2n}}
=\begin{cases}
\frac{1}{2^{\deg(Z)}}(-1)^{\langle\Lambda_1,\Lambda_2\rangle}, 
& \text{if $\Lambda_1,\Lambda_2\in\cals_Z$ for some special $Z$};\\
0, & \text{otherwise}.
\end{cases}
\]
Note that if $Z_\rmI=\binom{-}{-}$, i.e., if $Z=\binom{A}{B}$ with $A=B$,
then $\cals_Z=\{Z\}$ and $R_Z=\rho_Z$.
For this case, $Z$ is called \emph{degenerate} and we know that $Z^\rmt=Z$ and $\rho_Z=\sgn\cdot\rho_Z$.

Similar to the case of symplectic groups,
the set of unipotent almost characters of $\rmO^\pm_{2n}$ is denoted by
\[
\cala(\rmO^\pm_{2n})_1=\{\, R_\Lambda\mid\Lambda\in\cals_{\rmO^\pm_{2n}}\,\}.
\]

\begin{lem}\label{0208}
Let $Z$ be a special symbol of defect $0$.
The set $\cala_Z:=\{\,R_\Lambda\mid\Lambda\in\cals_Z\,\}$ forms
an orthonormal basis of $\calv_Z$.
\end{lem}
\begin{proof}
Let $M_1,M_2\subset Z_\rmI$.
Then by (\ref{0204}), we know that
\[
\langle R_{\Lambda_{M_1}},R_{\Lambda_{M_2}}\rangle
= \frac{1}{2^{2\deg(Z)}}\sum_{M\subset Z_\rmI}(-1)^{|M|(|M_1|+|M_2|)+|M\cap M_1|+|M\cap M_2|}
\]
If $M_1=M_2$, then clearly $\langle R_{\Lambda_{M_1}},R_{\Lambda_{M_2}}\rangle=1$.
Now suppose that $M_1\neq M_2$ and we have the following two situations:
\begin{enumerate}
\item Suppose that $|M_1|+|M_2|$ is even.
Let $a_1\in (M_2\smallsetminus M_1)\cup(M_1\smallsetminus M_2)$, which is non-empty.
It is not difficult to see that
\[
\{\, M\mid M\subset Z_\rmI\,\}
=\{\, M, M^*\mid M\subset Z_\rmI,\ a_1\not\in M\}
\]
where $M^*=M\cup\{a_1\}$.
Now
\[
|M^*\cap M_1|+|M^*\cap M_2|=|M\cap M_1|+|M\cap M_2|+1
\]
for any $M\subset Z_\rmI$ such that $a_1\not\in M$.
Hence, if $M_1\neq M_2$, we have
\begin{multline*}
\sum_{M\subset Z_\rmI}(-1)^{|M|(|M_1|+|M_2|)+|M\cap M_1|+|M\cap M_2|} \\
=\sum_{M\subset Z_\rmI,\ a_1\not\in M}\left[(-1)^{|M\cap M_1|+|M\cap M_2|}+(-1)^{|M^*\cap M_1|+|M^*\cap M_2|}\right]
=0.
\end{multline*}

\item Suppose that $|M_1|+|M_2|$ is odd.
Let $a_1\in (Z_\rmI\smallsetminus(M_1\cup M_2))\cup(M_1\cap M_2)$, which is non-empty.
It is not difficult to see that
\[
\{\, M\mid M\subset Z_\rmI\,\}
=\{\, M, M^*\mid M\subset Z_\rmI,\ a_1\not\in M\}
\]
where $M^*=M\cup\{a_1\}$.
Now
\[
|M^*\cap M_1|+|M^*\cap M_2|=\begin{cases}
|M\cap M_1|+|M\cap M_2|, & \text{if $a_1\in Z_\rmI\smallsetminus(M_1\cup M_2)$};\\
|M\cap M_1|+|M\cap M_2|+2, & \text{if $a_1\in M_1\cap M_2$}
\end{cases}
\]
for any $M\subset Z_\rmI$ such that $a_1\not\in M$.
Hence, if $M_1\neq M_2$, we have
\[
\sum_{M\subset Z_\rmI}(-1)^{|M|(|M_1|+|M_2|)+|M\cap M_1|+|M\cap M_2|}
=\sum_{M\subset Z_\rmI,\ a_1\not\in M}\left[(-1)^{|M|}+(-1)^{|M^*|}\right]
=0.
\]
\end{enumerate}
Therefore, the lemma is proved.
\end{proof}

\begin{exam}\label{0201}
It is known that $\cals_{\rmO^+_2}=\left\{\binom{1}{0},\binom{0}{1}\right\}$ and
$\cals_{\rmO^-_2}=\bigl\{\binom{-}{1,0},\binom{1,0}{-}\bigr\}$.
The symbol $Z=\binom{1}{0}$ is special of rank $1$ and degree $1$, and
$\cals_Z=\bigl\{\textstyle\binom{1}{0},\binom{0}{1},\binom{-}{1,0},\binom{1,0}{-}\bigr\}$.
Then the table of $(-1)^{\langle\Lambda_1,\Lambda_2\rangle}$ is given by
\[
\begin{tabular}{ccc|cccc}
\toprule
& & $\Lambda_{M_2}$ & $\binom{1}{0}$ & $\binom{0}{1}$ & $\binom{-}{1,0}$ & $\binom{1,0}{-}$ \\
& & $M_2$ & $\binom{-}{-}$ & $\binom{1}{0}$ & $\binom{0}{-}$ & $\binom{-}{0}$ \\
$\Lambda_{M_1}$ & $M_1$ & & & & \\
\midrule
$\binom{1}{0}$ & $\binom{-}{-}$ & & $\phantom{-}1$ & $\phantom{-}1$ & $\phantom{-}1$ & $\phantom{-}1$ \\
$\binom{0}{1}$ & $\binom{1}{0}$ & & $\phantom{-}1$ & $\phantom{-}1$ & $-1$ & $-1$ \\
$\binom{-}{1,0}$ & $\binom{0}{-}$ & & $\phantom{-}1$ & $-1$ & $\phantom{-}1$ & $-1$ \\
$\binom{1,0}{-}$ & $\binom{-}{1}$ & & $\phantom{-}1$ & $-1$ & $-1$ & $\phantom{-}1$ \\
\bottomrule
\end{tabular}
\]
Note that $\rho_{\binom{1}{0}}={\bf 1}_{\rmO^+_2}$,
$\rho_{\binom{0}{1}}=\sgn_{\rmO^+_2}$, $\rho_{\binom{-}{1,0}}={\bf 1}_{\rmO^-_2}$,
$\rho_{\binom{1,0}{-}}=\sgn_{\rmO^-_2}$.
Therefore we have the following four almost characters of $\rmO^\pm_2(q)$:
\begin{align*}
R_{\binom{1}{0}} &=\frac{1}{2}\left({\bf1}_{\rmO^+_2}+\sgn_{\rmO^+_2}+{\bf1}_{\rmO^-_2}+\sgn_{\rmO^-_2}\right), \\
R_{\binom{0}{1}} &=\frac{1}{2}\left({\bf1}_{\rmO^+_2}+\sgn_{\rmO^+_2}-{\bf1}_{\rmO^-_2}-\sgn_{\rmO^-_2}\right), \\
R_{\binom{-}{1,0}} &=\frac{1}{2}\left({\bf1}_{\rmO^+_2}-\sgn_{\rmO^+_2}+{\bf1}_{\rmO^-_2}-\sgn_{\rmO^-_2}\right), \\
R_{\binom{1,0}{-}} &=\frac{1}{2}\left({\bf1}_{\rmO^+_2}-\sgn_{\rmO^+_2}-{\bf1}_{\rmO^-_2}+\sgn_{\rmO^-_2}\right).
\end{align*}
\end{exam}

Let $\cals_{Z,0}$ be the set of symbols of defect $0$ in $\cals_Z$.
If $\Lambda\in\cals_{Z,0}$, then $\Upsilon(\Lambda)\in\calp_2(n)$.
Let $\chi_\Lambda$ be the irreducible character of $W_n$ parametrized by $\Upsilon(\Lambda)$.
Then it is known that
\[
R_\Lambda=\frac{1}{|W_n|}\left[\sum_{w\in W^+_n}\chi_\Lambda(w)R^{\rmO^+_{2n}}_{\bfT_w,1}
+\sum_{w\in W^-_n}\chi_\Lambda(w)R^{\rmO^-_{2n}}_{\bfT_w,1}\right]
\]
for $\Lambda\in\cals_{Z,0}$.
Here, $W^+_n$ is a subgroup of $W_n$ of index $2$, and $W^-_n=W_n\smallsetminus W^+_n$.
For $\epsilon=+$ or $-$, and $w\in W_n^\epsilon$,
it is known that $\bfT_w$ is a rational maximal torus in $\SO^\epsilon_{2n}$ and then
$R^{\rmO^\epsilon_{2n}}_{\bfT_w,1}:=\Ind^{\rmO^\epsilon_{2n}}_{\SO^\epsilon_{2n}}R^{\SO^\epsilon_{2n}}_{\bfT_w,1}$.

\section{Finite Howe Correspondence on Unipotent Almost Characters}
The purpose of this section is to prove Proposition~\ref{0301}.
All the techniques are adapted from \cite{pan-uniform} and so
some proofs might be somewhat sketchy.

\subsection{Finite Howe correspondence on unipotent characters}\label{0326}
Recall that for two partitions $[\lambda]=[\lambda_1,\lambda_2,\ldots,\lambda_m]$
(with $\lambda_1\geq\lambda_2\geq\cdots\geq\lambda_m$) and
$[\mu]=[\mu_1,\mu_2,\ldots,\mu_{m'}]$ (with $\mu_1\geq\mu_2\geq\cdots\geq\mu_{m'}$),
we denote $\lambda\preccurlyeq\mu$ if
\[
\mu_1\geq\lambda_1\geq\mu_2\geq\lambda_2\geq\cdots.
\]

Let $(\bfG,\bfG')=(\Sp_{2n},\rmO^\epsilon_{2n'})$ where $\epsilon=+$ or $-$.
For $\Lambda\in\cals_\bfG$ and $\Lambda'\in\cals_{\bfG'}$,
we write $\Upsilon(\Lambda)=\sqbinom{\lambda}{\mu}$ and $\Upsilon(\Lambda')=\sqbinom{\lambda'}{\mu'}$
and define
\begin{align*}
\calb_{\Sp_{2n},\rmO^+_{2n'}}
&=\bigl\{(\Lambda,\Lambda')\in\cals_{\Sp_{2n}}\times\cals_{\rmO^+_{2n'}}\mid\mu\preccurlyeq\lambda',\ \mu'\preccurlyeq\lambda,\
{\rm def}(\Lambda')=-{\rm def}(\Lambda)+1\,\bigr\} \\
\calb_{\Sp_{2n},\rmO^-_{2n'}}
&=\bigl\{(\Lambda,\Lambda')\in\cals_{\Sp_{2n}}\times\cals_{\rmO^-_{2n'}}\mid\lambda'\preccurlyeq\mu,\ \lambda\preccurlyeq\mu',\
{\rm def}(\Lambda')=-{\rm def}(\Lambda)-1\,\bigr\} \\
\calb_{\Sp_{2n},\rmO^\pm_{2n'}}
&=\calb_{\Sp_{2n},\rmO^+_{2n'}}\cup\calb_{\Sp_{2n},\rmO^-_{2n'}}\subset\cals_{\Sp_{2n}}\times\cals_{\rmO^\pm_{2n'}}.
\end{align*}

The following proposition is proved in \cite{pan-finite-unipotent}:
\begin{prop}
Let $n,n'$ be non-negative integers.
Then we have
\[
\omega_{\Sp_{2n},\rmO^\epsilon_{2n'},1}
=\sum_{(\Lambda,\Lambda')\in\calb_{\Sp_{2n},\rmO^\epsilon_{2n'}}}\rho_\Lambda\otimes \rho_{\Lambda'}
\]
for $\epsilon=+$ or $-1$.
\end{prop}

From the proposition, we conclude that
\begin{equation}\label{0325}
\omega_{\Sp_{2n},\rmO^\pm_{2n'},1}
=\sum_{(\Lambda,\Lambda')\in\calb_{\Sp_{2n},\rmO^\pm_{2n'}}}\rho_\Lambda\otimes \rho_{\Lambda'}.
\end{equation}

\begin{exam}\label{0322}
Consider the pair $(\Sp_4,\rmO^\pm_2)$.
It is not difficult to see that the table for $\calb_{\Sp_4,\rmO^\pm_2}$ is
\[
\begin{tabular}{cc|c|c|ccccc}
\toprule
& $\cals_{\Sp_4}$ & $\binom{2}{-}$ & $\binom{2,1,0}{2,1}$ & $\binom{2,0}{1}$ & $\binom{2,1}{0}$ & $\binom{1,0}{2}$ & $\binom{-}{2,1,0}$ \\
$\cals_{\rmO^\pm_2}$ & & $\sqbinom{2}{-}$ & $\sqbinom{-}{1,1}$ & $\sqbinom{1}{1}$ & $\sqbinom{1,1}{-}$ & $\sqbinom{-}{2}$ & $\sqbinom{-}{-}$ \\
\midrule
$\binom{1}{0}$ & $\sqbinom{1}{-}$ & \checkmark & & \checkmark & & & \\
$\binom{0}{1}$ & $\sqbinom{-}{1}$ & \checkmark & & & \checkmark & & \\
\midrule
$\binom{-}{1,0}$ & $\sqbinom{-}{-}$ & & & & & \checkmark & \\
$\binom{1,0}{-}$ & $\sqbinom{-}{-}$ & & & & & & \checkmark  \\
\bottomrule
\end{tabular}
\]
This means that
\begin{align*}
\omega_{\Sp_4,\rmO^+_2,1}
&=\rho_{\binom{2}{-}}\otimes\bigl[\rho_{\binom{1}{0}}+\rho_{\binom{0}{1}}\bigr]
+\rho_{\binom{2,0}{1}}\otimes\rho_{\binom{1}{0}}+\rho_{\binom{2,1}{0}}\otimes\rho_{\binom{0}{1}}, \\
\omega_{\Sp_4,\rmO^-_2,1}
&=\rho_{\binom{-}{2,1,0}}\otimes\rho_{\binom{1,0}{-}}+\rho_{\binom{1,0}{2}}\otimes\rho_{\binom{-}{1,0}}.
\end{align*}
From Example~\ref{0201}, for $\rmO^\pm_2$, we know that
\[
\rho_{\binom{1}{0}}+\rho_{\binom{0}{1}}=R_{\binom{1}{0}}+R_{\binom{0}{1}}.
\]
For $\Sp_4$, it is known that
\begin{align*}
\rho_{\binom{2}{-}}
&= R_{\binom{2}{-}}, \\
\rho_{\binom{2,0}{1}}
&=\frac{1}{2}\left[R_{\binom{2,1}{0}}+R_{\binom{2,0}{1}}+R_{\binom{1,0}{2}}+R_{\binom{-}{2,1,0}}\right], \\
\rho_{\binom{2,1}{0}}
&=\frac{1}{2}\left[R_{\binom{2,1}{0}}+R_{\binom{2,0}{1}}-R_{\binom{1,0}{2}}-R_{\binom{-}{2,1,0}}\right].
\end{align*}
It is not difficult to check that
\begin{multline*}
\rho_{\binom{2,0}{1}}\otimes\rho_{\binom{1}{0}}+\rho_{\binom{2,1}{0}}\otimes\rho_{\binom{0}{1}}
+\rho_{\binom{-}{2,1,0}}\otimes\rho_{\binom{1,0}{-}}+\rho_{\binom{1,0}{2}}\otimes\rho_{\binom{-}{1,0}} \\
=R_{\binom{2,0}{1}}\otimes R_{\binom{1}{0}}+R_{\binom{2,1}{0}}\otimes R_{\binom{0}{1}}
+R_{\binom{-}{2,1,0}}\otimes R_{\binom{1,0}{-}}+R_{\binom{1,0}{2}}\otimes R_{\binom{-}{1,0}}.
\end{multline*}
Therefore we conclude that
\[
\omega_{\Sp_4,\rmO^\pm_2,1}=\sum_{(\Lambda,\Lambda')\in\calb_{\Sp_4,\rmO^\pm_2}}\rho_\Lambda\otimes\rho_{\Lambda'}
=\sum_{(\Lambda,\Lambda')\in\calb_{\Sp_4,\rmO^\pm_2}}R_\Lambda\otimes R_{\Lambda'}.
\]
If we consider two dual pairs $(\Sp_4,\rmO^+_2)$ and $(\Sp_4,\rmO^-_2)$ separately,
it is not difficult to see that
\begin{align*}
\omega_{\Sp_4,\rmO^+_2,1}
&=\sum_{(\Lambda,\Lambda')\in\calb_{\Sp_4,\rmO^+_2}}\rho_\Lambda\otimes\rho_{\Lambda'}
\neq\sum_{(\Lambda,\Lambda')\in\calb_{\Sp_4,\rmO^+_2}}R_\Lambda\otimes R_{\Lambda'}; \\
\omega_{\Sp_4,\rmO^-_2,1}
&=\sum_{(\Lambda,\Lambda')\in\calb_{\Sp_4,\rmO^-_2}}\rho_\Lambda\otimes\rho_{\Lambda'}
\neq\sum_{(\Lambda,\Lambda')\in\calb_{\Sp_4,\rmO^-_2}}R_\Lambda\otimes R_{\Lambda'}.
\end{align*}
This means Theorem~\ref{0101} does not hold if we consider two dual pairs $(\Sp_{2n},\rmO^+_{2n'})$ and
$(\Sp_{2n},\rmO^-_{2n'})$ separately.
\end{exam}

For special symbols $Z,Z'$ of ranks $n,n'$ and of defects $1,0$ respectively, we define
\begin{align*}
\calb^\epsilon_{Z,Z'}
&= \calb_{\Sp_{2n},\rmO^\epsilon_{2n'}}\cap(\cals_Z\times\cals^\epsilon_{Z'}), \\
\calb_{Z,Z'}=\calb^+_{Z,Z'}\cup\calb^-_{Z,Z'}
&= \calb_{\Sp_{2n},\rmO^\pm_{2n'}}\cap(\cals_Z\times\cals_{Z'}).
\end{align*}
Then we have
\begin{align*}
\calb_{\Sp_{2n},\rmO^\pm_{2n'}}
&=\bigcup_{Z,Z'}\calb_{Z,Z'}, \\
\sum_{(\Lambda,\Lambda')\in\calb_{\Sp_{2n},\rmO^\pm_{2n'}}}\rho_\Lambda\otimes\rho_{\Lambda'}
&=\sum_{Z,Z'}\;\sum_{(\Lambda,\Lambda')\in\calb_{Z,Z'}}\rho_\Lambda\otimes\rho_{\Lambda'}, \\
\sum_{(\Lambda,\Lambda')\in\calb_{\Sp_{2n},\rmO^\pm_{2n'}}} R_\Lambda\otimes R_{\Lambda'}
&=\sum_{Z,Z'}\;\sum_{(\Lambda,\Lambda')\in\calb_{Z,Z'}} R_\Lambda\otimes R_{\Lambda'}
\end{align*}
where $Z,Z'$ run over special symbols of ranks $n,n'$ and of defects $1,0$ respectively.

\begin{prop}\label{0301}
Consider the dual pair $(\bfG,\bfG')=(\Sp_{2n},\rmO^\pm_{2n'})$,
and let $Z,Z'$ be special symbols of ranks $n,n'$ and defects $1,0$ respectively.
Then
\[
\sum_{(\Lambda,\Lambda')\in\calb_{Z,Z'}}\rho_\Lambda\otimes\rho_{\Lambda'}
=\sum_{(\Lambda,\Lambda')\in\calb_{Z,Z'}}R_\Lambda\otimes R_{\Lambda'}.
\]
\end{prop}
From Subsection~\ref{0327} to the end of this section,
we are going to prove Proposition~\ref{0301}.
We shall consider some special cases first
and finally reduce the general case to special cases in Subsection~\ref{0310}.

\begin{cor}\label{0306}
Consider the dual pair $(\bfG,\bfG')=(\Sp_{2n},\rmO^\pm_{2n'})$.
Then
\[
\omega_{\Sp_{2n},\rmO^\pm_{2n'},1}
=\sum_{(\Lambda,\Lambda')\in\calb_{\Sp_{2n},\rmO^\pm_{2n'}}}\rho_\Lambda\otimes\rho_{\Lambda'}
=\sum_{(\Lambda,\Lambda')\in\calb_{\Sp_{2n},\rmO^\pm_{2n'}}}R_\Lambda\otimes R_{\Lambda'}.
\]
\end{cor}
\begin{proof}
By (\ref{0325}) and Proposition~\ref{0301}, we have
\begin{align*}
\omega_{\Sp_{2n},\rmO^\pm_{2n'},1}
=\sum_{(\Lambda,\Lambda')\in\calb_{\Sp_{2n},\rmO^\pm_{2n'}}}\rho_\Lambda\otimes\rho_{\Lambda'}
&=\sum_{Z,Z'}\;\sum_{(\Lambda,\Lambda')\in\calb_{Z,Z'}}\rho_\Lambda\otimes\rho_{\Lambda'} \\
&=\sum_{Z,Z'}\;\sum_{(\Lambda,\Lambda')\in\calb_{Z,Z'}}R_\Lambda\otimes R_{\Lambda'} \\
&=\sum_{(\Lambda,\Lambda')\in\calb_{\Sp_{2n},\rmO^\pm_{2n'}}}R_\Lambda\otimes R_{\Lambda'}
\end{align*}
where $Z,Z'$ run over special symbols of ranks $n,n'$ and of defects $1,0$ respectively.
\end{proof}

From Corollary~\ref{0306},
we have a relation between the set of unipotent almost characters of $\Sp_{2n}(q)$ and the set
of unipotent almost characters of $\rmO^\pm_{2n'}(q)$:
\[
\widetilde\Theta_{\Sp_{2n},\rmO^\pm_{2n'},1}
=\bigl\{\,(R_\Lambda,R_{\Lambda'})\in\cala(\Sp_{2n})_1\times\cala(\rmO^\pm_{2n'})_1
\mid(\Lambda,\Lambda')\in\calb_{\Sp_{2n},\rmO^\pm_{2n'}}\,\bigr\}.
\]

\subsection{Basic case I}\label{0327}
For a non-negative integer $m$, we define two special symbols
\begin{equation}
Z_{(m)}=\binom{2m,2m-2,\ldots,0}{2m-1,2m-3,\ldots,1}, \qquad
Z'_{(m)}=\binom{2m-1,2m-3,\ldots,1}{2m-2,2m-4,\ldots,0}
\end{equation}
of ranks $m(m+1),m^2$ and defects $1,0$ respectively.
In this subsection, we always assume that $Z=Z_{(m)}$ and $Z'=Z'_{(m)}$,
and we want to confirm Proposition~\ref{0301} for this special case.
Now $\bfG=\Sp_{2m(m+1)}$, $\bfG'=\rmO^\pm_{2m^2}$, $\deg(Z_{(m)})=\deg({Z'_{(m)}})=m$ and
$|\cals_Z|=|\cals_{Z'}|=2^{2m}$.
Note that $Z_\rmI=Z$ and $Z'_\rmI=Z'$.
Define a mapping $\theta\colon\cals_{Z'}\rightarrow\cals_Z$ given by
\begin{equation}\label{0323}
\theta(\Lambda')=\begin{cases}
\binom{2m}{-}\cup\Lambda'^\rmt, & \text{if ${\rm def}(\Lambda')\equiv 0\pmod 4$};\\
\binom{-}{2m}\cup\Lambda'^\rmt, & \text{if ${\rm def}(\Lambda')\equiv 2\pmod 4$}.
\end{cases}
\end{equation}
It is not difficult to see that ${\rm def}(\theta(\Lambda'))\equiv 1\pmod 4$ and $\theta$ is a bijection.

\begin{lem}\label{0304}
Let $Z=Z_{(m)}$ and $Z'=Z'_{(m)}$.
Then we have
$\calb_{Z,Z'}=\{\,(\theta(\Lambda'),\Lambda')\mid\Lambda'\in\cals_{Z'}\,\}$.
\end{lem}
\begin{proof}
From \cite{pan-uniform} example~2.27, it is known that
$\calb^\epsilon_{Z,Z'}=\{\,(\theta(\Lambda'),\Lambda')\mid\Lambda'\in\cals^\epsilon_{Z'}\,\}$
for $\epsilon=+$ or $-$.
Now $\cals_{Z'}=\cals^+_{Z'}\cup\cals^-_{Z'}$ and $\calb_{Z,Z'}=\calb^+_{Z,Z'}\cup\calb^-_{Z,Z'}$.
So the lemma is proved.
\end{proof}

\begin{lem}\label{0305}
Let $Z=Z_{(m)}$ and $Z'=Z'_{(m)}$.
Then we have
$\langle\Lambda'_1,\Lambda'_2\rangle=\langle\theta(\Lambda'_1),\theta(\Lambda'_2)\rangle$
for any $\Lambda'_1,\Lambda'_2\in\cals_{Z'}$.
\end{lem}
\begin{proof}
Write $\Lambda'_1=\Lambda_{M'_1}$ and $\Lambda'_2=\Lambda_{M'_2}$ for some $M'_1,M'_2\subset Z'_\rmI$ where
$\Lambda_{M'_i}$ is given in (\ref{0202}).
Note that $M_1'^\rmt,M_2'^\rmt\subset Z_\rmI$.
Then
\[
\theta(\Lambda_{M'_i})=\begin{cases}
\Lambda_{M_i'^\rmt}, & \text{if $|M'_i|$ is even};\\
\Lambda_{(M_i'^*)^\rmt}, & \text{if $|M'_i|$ is odd}
\end{cases}
\]
where $M_i'^*=\binom{-}{2m}\cup M'_i$.
Note that $|M'^*_i\cap M'_j|=|M'_i\cap M'_j|$ and $|M'^*_i\cap M'^*_j|=|M'_i\cap M'_j|+1$.
Now we have the following three cases:
\begin{enumerate}
\item If both $|M'_1|,|M'_2|$ are even, then by (\ref{0203}) and (\ref{0204})
we have
\begin{align*}
\langle\theta(\Lambda_{M'_1}),\theta(\Lambda_{M'_2})\rangle
=\langle\Lambda_{M_1'^\rmt},\Lambda_{M_2'^\rmt}\rangle
&\equiv |M_1'^\rmt\cap M_2'^\rmt|\pmod 2  \\
&\equiv |M'_1||M'_2|+|M'_1\cap M'_2|\pmod 2 \\
&=\langle\Lambda_{M'_1},\Lambda_{M'_2}\rangle.
\end{align*}

\item If $|M'_1|+|M'_2|$ is odd, then exactly one of $|M_1'|,|M_2'|$ is even and
\begin{align*}
\langle\theta(\Lambda_{M'_1}),\theta(\Lambda_{M'_2})\rangle
=\langle\Lambda_{(M_1'^*)^\rmt},\Lambda_{M_2'^\rmt}\rangle
&\equiv |(M_1'^*)^\rmt\cap M_2'^\rmt|\pmod 2  \\
&\equiv |M_1'^*||M'_2|+|M'_1\cap M'_2|\pmod 2 \\
&=\langle\Lambda_{M'_1},\Lambda_{M'_2}\rangle.
\end{align*}

\item If both $|M'_1|,|M'_2|$ are odd, then
\begin{align*}
\langle\theta(\Lambda_{M'_1}),\theta(\Lambda_{M'_2})\rangle
=\langle\Lambda_{M_1'^\rmt},\Lambda_{M_2'^\rmt}\rangle
&\equiv |(M_1'^*)^\rmt\cap(M_2'^*)^\rmt|\pmod 2  \\
&\equiv |M'_1||M'_2|+1+|M'_1\cap M'_2|\pmod 2 \\
&=\langle\Lambda_{M'_1},\Lambda_{M'_2}\rangle.
\end{align*}
\end{enumerate}
\end{proof}

\begin{lem}\label{0302}
Let $Z=Z_{(m)}$ and $Z'=Z'_{(m)}$.
Then we have
\[
\sum_{(\Lambda,\Lambda')\in\calb_{Z,Z'}}\rho_\Lambda\otimes\rho_{\Lambda'}
=\sum_{(\Lambda,\Lambda')\in\calb_{Z,Z'}}R_\Lambda\otimes R_{\Lambda'}.
\]
\end{lem}
\begin{proof}
Write
\[
\sum_{(\Lambda,\Lambda')\in\calb_{Z,Z'}}\rho_\Lambda\otimes\rho_{\Lambda'}
=\sum_{\Lambda_1\in\cals_Z,\ \Lambda'_1\in\cals_{Z'}}m_{\Lambda_1,\Lambda_1'}R_{\Lambda_1}\otimes R_{\Lambda'_1}
\]
where $m_{\Lambda_1,\Lambda'_1}\in\bbC$.
Then by (\ref{0205}), (\ref{0206}), Lemma~\ref{0304}, Lemma~\ref{0305} and the proof of Lemma~\ref{0207},
we have
\begin{align*}
m_{\Lambda_1,\Lambda_1'}
=\left\langle\sum_{\Lambda'\in\cals_{Z'}}\rho_{\theta(\Lambda')}\otimes\rho_{\Lambda'},
R_{\Lambda_1}\otimes R_{\Lambda'_1}\right\rangle_{\bfG\times\bfG'}
&=\frac{1}{2^{2m}}\sum_{\Lambda'\in\cals_{Z'}}(-1)^{\langle\theta(\Lambda'),\Lambda_1\rangle+\langle\Lambda',\Lambda'_1\rangle} \\
&=\frac{1}{2^{2m}}\sum_{\Lambda'\in\cals_{Z'}}(-1)^{\langle\Lambda',\theta^{-1}(\Lambda_1)\rangle+\langle\Lambda',\Lambda'_1\rangle} \\
&=\begin{cases}
1, & \text{if $\Lambda_1=\theta(\Lambda'_1)$};\\
0, & \text{if $\Lambda_1\neq\theta(\Lambda'_1)$}.
\end{cases}
\end{align*}
Now $\Lambda_1=\theta(\Lambda'_1)$ if and only if $(\Lambda_1,\Lambda_1')\in\calb_{Z,Z'}$.
Therefore we conclude that $m_{\Lambda_1,\Lambda_1'}=1$ if $(\Lambda_1,\Lambda_1')\in\calb_{Z,Z'}$
and $m_{\Lambda_1,\Lambda_1'}=0$ if $(\Lambda_1,\Lambda_1')\not\in\calb_{Z,Z'}$.
\end{proof}

\subsection{Basic case II}
In this subsection we assume that $Z=Z_{(m)}$ and $Z'=Z'_{(m+1)}$ for some non-negative integer $m$,
and we want to prove Proposition~\ref{0301} for this special case by reducing it to the case that
considered in the previous subsection.

Now $\bfG=\Sp_{2m(m+1)}$, $\bfG'=\rmO^\pm_{2(m+1)^2}$,
and we have $|\cals_Z|=2^{2m}$ and $|\cals_{Z'}|=2^{2(m+1)}$.
We define two injective mappings $\theta^+,\theta^-\colon\cals_Z\rightarrow\cals_{Z'}$ given by
\begin{equation}\label{0312}
\theta^\epsilon(\Lambda)=\begin{cases}
\binom{2m+1}{-}\cup\Lambda^\rmt, & \text{if $\epsilon=+$}; \\
\binom{-}{2m+1}\cup\Lambda^\rmt, & \text{if $\epsilon=-$}.
\end{cases}
\end{equation}
It is clear that the image of $\theta^\epsilon$ is in $\cals^\epsilon_{Z'}$, in fact
\begin{align*}
\theta^+(\cals_Z) &=\{\,\Lambda_{M'}\mid M'\subset Z'\smallsetminus\Psi',\ |M'|\text{ even}\,\}\subset\cals^+_{Z'} \\
\theta^-(\cals_Z) &=\{\,\Lambda_{M'\cup\Psi'}\mid M'\subset Z'\smallsetminus\Psi',\ |M'|\text{ even}\,\}\subset\cals^-_{Z'}
\end{align*}
where $\Psi'=\binom{2m+1}{-}$.
Moreover, from \cite{pan-uniform} example~2.26, we know that
\begin{equation}\label{0308}
\calb_{Z,Z'}=\{\,(\Lambda,\theta^+(\Lambda)),(\Lambda,\theta^-(\Lambda))\mid\Lambda\in\cals_Z\,\}
\end{equation}

\begin{lem}
Let $Z=Z_{(m)}$ and $Z'=Z'_{(m+1)}$.
Then we have
\[
\langle\Lambda_1,\Lambda_2\rangle
=\langle\theta^+(\Lambda_1),\theta^+(\Lambda_2)\rangle
=\langle\theta^-(\Lambda_1),\theta^-(\Lambda_2)\rangle
\]
for any $\Lambda_1,\Lambda_2\in\cals_Z$.
\end{lem}
\begin{proof}
Let $M_1,M_2\subset Z_\rmI$ such that both $|M_1|,|M_2|$ are even.
Now by (\ref{0204}), we have
\begin{align*}
\langle\theta^+(\Lambda_{M_1}),\theta^+(\Lambda_{M_2})\rangle
&=|M_1^\rmt||M_2^\rmt|+|M_1^\rmt\cap M_2^\rmt|\pmod 2 \\
&=|M_1\cap M_2|\pmod 2
=\langle\Lambda_{M_1},\Lambda_{M_2}\rangle.
\end{align*}
Similarly,
\begin{align*}
\langle\theta^-(\Lambda_{M_1}),\theta^-(\Lambda_{M_2})\rangle
&=|M_1^\rmt\cup\Psi'||M_2^\rmt\cup\Psi'|+|(M_1^\rmt\cup\Psi')\cap(M_2^\rmt\cup\Psi')|\pmod 2 \\
&=|M_1\cap M_2|\pmod 2
=\langle\Lambda_{M_1},\Lambda_{M_2}\rangle
\end{align*}
where $\Psi'=\binom{2m+1}{-}$.
\end{proof}

Let $\bfG'^{(1)}=\rmO^\pm_{2m^2}$, $Z'^{(1)}=Z'_{(m)}\in\cals_{\bfG'^{(1)}}$, and $\Psi'=\binom{2m+1}{-}$.
Define a mapping
\begin{align*}
f'\colon\{\,M',M'\cup\Psi'\mid M'\subset Z'\smallsetminus\Psi',\ |M'|\ \text{even}\} &\rightarrow\{\,M''\mid M''\subset Z'^{(1)}\,\} \\
M',M'\cup\Psi' &\mapsto M'\smallsetminus\textstyle\binom{-}{2m}.
\end{align*}
It is not difficult to see that $f'$ is a surjective two-to-one mapping.

\begin{lem}
Let $Z=Z_{(m)}$, $Z'=Z'_{(m+1)}$ and $Z'^{(1)}=Z'_{(m)}$.
Then $(\Lambda_M,\Lambda_{M'})\in\calb_{Z,Z'}$ if and only if\/
$(\Lambda_M,\Lambda_{f'(M')})\in\calb_{Z,Z'^{(1)}}$.
\end{lem}
\begin{proof}
From (\ref{0308}) and Lemma~\ref{0304}, we know that
\begin{align*}
\calb_{Z,Z'} &=\{\,(\Lambda_M,\Lambda_{M^\rmt}),(\Lambda_M,\Lambda_{M^\rmt\cup\Psi'})\mid M\subset Z,\ |M|\ \text{even}\} \\
\calb_{Z,Z'^{(1)}} &=\{\,(\theta(\Lambda_{M''}),\Lambda_{M''})\mid M''\subset Z'^{(1)}\,\}
\end{align*}
where $\theta$ is given in (\ref{0323}).
Note that for $M\subset Z$ with $|M|$ even,
we have $M^\rmt\subset Z'\smallsetminus\Psi'$ with $|M^\rmt|$ even.
So $f'(M^\rmt)=f'(M^\rmt\cup\Psi')=M^\rmt\smallsetminus\binom{-}{2m}$.
Now from (\ref{0312}), we see that
\begin{itemize}
\item if $2m\not\in M$, then $f'(M^\rmt)=f'(M^\rmt\cup\Psi')=M^\rmt$ and $\theta(\Lambda_{M^\rmt})=\Lambda_M$;

\item if $2m\in M$, then $f'(M^\rmt)=f'(M^\rmt\cup\Psi')=M^\rmt\smallsetminus\binom{-}{2m}$ and
$\theta\bigl(\Lambda_{M^\rmt\smallsetminus\binom{-}{2m}}\bigr)=\Lambda_M$.
\end{itemize}
Therefore, the lemma is proved.
\end{proof}

\begin{exam}
Let $Z=Z_{(1)}=\binom{2,0}{1}\in\cals_{\Sp_4}$, $Z'=Z'_{(2)}=\binom{3,1}{2,0}\in\cals_{\rmO^\pm_8}$.
Then the relation $\calb_{Z,Z'}$ is given by the table:
\[
\begin{tabular}{ccc|cccc|cccc}
\toprule
& $\cals_{\rmO^\pm_8}$ & $\Lambda_{M'}$ & $\binom{3,1}{2,0}$ & $\binom{3,0}{2,1}$ & $\binom{3,2}{1,0}$ & $\binom{3,2,1,0}{-}$ & $\binom{1}{3,2,0}$ & $\binom{0}{3,2,1}$ & $\binom{2}{3,1,0}$ & $\binom{2,1,0}{3}$ \\
$\cals_{\Sp_4}$ & & $M'$ & $\binom{-}{-}$ & $\binom{1}{0}$ & $\binom{1}{2}$ & $\binom{-}{2,0}$ & $\binom{3}{-}$ & $\binom{3,1}{0}$
& $\binom{3,1}{2}$ & $\binom{3}{2,0}$ \\
$\Lambda_M$ & $M$ & & & & & & & \\
\midrule
$\binom{2,0}{1}$ & $\binom{-}{-}$ & & \checkmark & & & & \checkmark & & & \\
$\binom{2,1}{0}$ & $\binom{0}{1}$ & & & \checkmark & & & & \checkmark & & \\
$\binom{1,0}{2}$ & $\binom{2}{1}$ & & & & \checkmark & & & & \checkmark & \\
$\binom{-}{2,1,0}$ & $\binom{2,0}{-}$ & & & & & \checkmark & & & & \checkmark \\
\bottomrule
\end{tabular}
\]
Now $Z'^{(1)}=\binom{1}{0}\in\cals_{\rmO^\pm_2}$ and
$f'(M')=M'\smallsetminus\binom{-}{2}$ for any $M'\subset\binom{1}{2,0}$ with $|M'|$ even.
Then the relation $\calb_{Z,Z'^{(1)}}$ is given by the table:
\[
\begin{tabular}{ccc|cc|cc}
\toprule
& $\cals_{\rmO^\pm_2}$ & $\Lambda_{f'(M')}$ & $\binom{1}{0}$ & $\binom{0}{1}$ & $\binom{-}{1,0}$ & $\binom{1,0}{-}$  \\
$\cals_{\Sp_4}$ & & $f'(M')$ & $\binom{-}{-}$ & $\binom{1}{0}$ & $\binom{1}{-}$ & $\binom{-}{0}$  \\
$\Lambda_M$ & $M$ & & & \\
\midrule
$\binom{2,0}{1}$ & $\binom{-}{-}$ & & \checkmark & & & \\
$\binom{2,1}{0}$ & $\binom{0}{1}$ & & & \checkmark & & \\
$\binom{1,0}{2}$ & $\binom{2}{1}$ & & & & \checkmark & \\
$\binom{-}{2,1,0}$ & $\binom{2,0}{-}$ & & & & & \checkmark \\
\bottomrule
\end{tabular}
\]
\end{exam}

For $M'\subset Z'\smallsetminus\Psi'$ with $|M'|$ even where $\Psi'=\binom{2m+1}{-}$,
we define
\[
\rho^{(1)}_{\Lambda_{M'}}=\frac{1}{\sqrt 2}(\rho_{\Lambda_{M'}}+\rho_{\Lambda_{M'\cup\Psi'}})\quad\text{and}\quad
R^{(1)}_{\Lambda_{M'}}=\frac{1}{\sqrt 2}(R_{\Lambda_{M'}}+R_{\Lambda_{M'\cup\Psi'}}).
\]
Let $\calv^{(1)}_{Z'}\subset\calv_{Z'}$ be the subspace spanned by
\[
\left\{\,\rho^{(1)}_{\Lambda_{M'}}\mid M'\subset Z'\smallsetminus\Psi',\ |M'|\ \text{even}\,\right\}.
\]
Then $f'$ induces a linear transformation
\[
\widetilde f'\colon\calv^{(1)}_{Z'}\rightarrow\calv_{Z'^{(1)}}\quad\text{by}\quad
\rho^{(1)}_{\Lambda_{M'}}\mapsto\rho_{\Lambda_{f'(M')}}.
\]

\begin{lem}\label{0309}
For $M_1',M_2'\subset Z'\smallsetminus\Psi'$ with $|M'_1|,|M'_2|$ both even,
we have
\[
\Bigl\langle \rho^{(1)}_{\Lambda_{M_1'}},R^{(1)}_{\Lambda_{M_2'}}\Bigr\rangle_{\bfG'}
=\Bigl\langle\rho_{\Lambda_{f'(M_1')}},R_{\Lambda_{f'(M_2')}}\Bigr\rangle_{\bfG'^{(1)}}.
\]
\end{lem}
\begin{proof}
Note that both $|M'_1|,|M'_2|$ are even, and we have
\begin{align*}
|M_1'||M'_2| &\equiv |M_1'\cup\Psi'||M'_2|\equiv |M_1'||M'_2\cup\Psi'|\equiv |M_1'\cup\Psi'||M'_2\cup\Psi'|+1\equiv 0\pmod 2, \\
|M_1'\cap M'_2| &= |(M_1'\cup\Psi')\cap M'_2|=|M_1'\cap(M'_2\cup\Psi')|=|(M_1'\cup\Psi')\cap(M'_2\cup\Psi')|-1.
\end{align*}
Therefore
\begin{align*}
& \Bigl\langle \rho^{(1)}_{\Lambda_{M_1'}},R^{(1)}_{\Lambda_{M_2'}}\Bigr\rangle_{\bfG'} \\
&= \frac{1}{2}\left[\bigl\langle \rho_{\Lambda_{M_1'}},R_{\Lambda_{M_2'}}\bigr\rangle_{\bfG'}
+\bigl\langle \rho_{\Lambda_{M_1'\cup\Psi'}},R_{\Lambda_{M_2'}}\bigr\rangle_{\bfG'}
+\bigl\langle \rho_{\Lambda_{M_1'}},R_{\Lambda_{M_2'\cup\Psi'}}\bigr\rangle_{\bfG'}
+\bigl\langle \rho_{\Lambda_{M_1'\cup\Psi'}},R_{\Lambda_{M_2'\cup\Psi'}}\bigr\rangle_{\bfG'}\right] \\
&= 2\,\bigl\langle \rho_{\Lambda_{M_1'}},R_{\Lambda_{M_2'}}\bigr\rangle_{\bfG'}.
\end{align*}
Moreover,
\begin{itemize}
\item if $2m\not\in M'_1\cap M'_2$,
then
\[
|f'(M'_1)||f'(M'_2)|\equiv |M'_1||M'_2|\equiv 0\pmod 2\quad\text{and}\quad |f'(M'_1)\cap f'(M'_2)|=|M'_1\cap M'_2|;
\]

\item if $2m\in M'_1\cap M'_2$,
then
\[
|f'(M'_1)||f'(M'_2)|\equiv |M'_1||M'_2|+1\pmod2\quad\text{and}\quad |f'(M'_1)\cap f'(M'_2)|=|M'_1\cap M'_2|-1.
\]
\end{itemize}
Now $\deg(Z')=\deg(Z'^{(1)})+1$, we have
\begin{align*}
\Bigl\langle\rho_{\Lambda_{f'(M_1')}},R_{\Lambda_{f'(M_2')}}\Bigr\rangle_{\bfG'^{(1)}}
&=2\,\bigl\langle \rho_{\Lambda_{M_1'}},R_{\Lambda_{M_2'}}\bigr\rangle_{\bfG'},
\end{align*}
and so the lemma is proved.
\end{proof}

\begin{cor}
The linear transformation $\widetilde f'\colon\calv^{(1)}_{Z'}\rightarrow\calv_{Z'^{(1)}}$
is an isometry, $R^{(1)}_{\Lambda_{M'}}\in\calv^{(1)}_{Z'}$, and
$\widetilde f'(R^{(1)}_{\Lambda_{M'}})=R_{\Lambda_{f'(M')}}$
for any $M'\subset Z'\smallsetminus\Psi'$ with $|M'|$ even.
\end{cor}
\begin{proof}
The linear transformation $\widetilde f'$ maps the orthonormal basis
\[
\left\{\,\rho^{(1)}_{\Lambda_{M'}}\mid M'\subset Z'\smallsetminus\Psi',\ |M'|\ \text{even}\,\right\}
\]
for $\calv^{(1)}_{Z'}$ bijectively
to the orthonormal basis $\{\,\rho_{\Lambda_{f'(M')}}\mid M'\subset Z'^{(1)}\,\}$ for $\calv_{Z'^{(1)}}$,
so $\widetilde f'$ is an isometry.
From the proof of Lemma~\ref{0309}, we see that
\[
\bigl\langle \rho_{\Lambda_{M_1'}},R_{\Lambda_{M_2'}}\bigr\rangle_{\bfG'}
=\bigl\langle \rho_{\Lambda_{M_1'\cup\Psi'}},R_{\Lambda_{M_2'\cup\Psi'}}\bigr\rangle_{\bfG'}
\]
for any $M'_1,M'_2\subset Z'\smallsetminus\Psi'$ with $|M'_1|,|M'_2|$ both even.
This implies that $R^{(1)}_{\Lambda_{M'}}\in\calv^{(1)}_{Z'}$.
Moreover, we have $\widetilde f'(R^{(1)}_{\Lambda_{M'}})=R_{\Lambda_{f'(M')}}$ by Lemma~\ref{0309}.
\end{proof}

\begin{lem}\label{0307}
Let $m$ be a non-negative integer,
and let $Z=Z_{(m)}$ and $Z'=Z'_{(m+1)}$.
Then
\[
\sum_{(\Lambda,\Lambda')\in\calb_{Z,Z'}}\rho_\Lambda\otimes\rho_{\Lambda'}
=\sum_{(\Lambda,\Lambda')\in\calb_{Z,Z'}}R_\Lambda\otimes R_{\Lambda'}.
\]
\end{lem}
\begin{proof}
By (\ref{0308}), we have
\begin{align*}
\sum_{(\Lambda,\Lambda')\in\calb_{Z,Z'}}\rho_\Lambda\otimes\rho_{\Lambda'}
=\sum_{M\subset Z_\rmI,\ |M|\text{ even}}\rho_{\Lambda_M}\otimes(\rho_{\Lambda_{M^\rmt}}+\rho_{\Lambda_{M^\rmt\cup\Psi'}})
=\sqrt 2\sum_{M\subset Z_\rmI,\ |M|\text{ even}}\rho_{\Lambda_M}\otimes\rho^{(1)}_{\Lambda_{M^\rmt}} \\
\sum_{(\Lambda,\Lambda')\in\calb_{Z,Z'}}R_\Lambda\otimes R_{\Lambda'}
=\sum_{M\subset Z_\rmI,\ |M|\text{ even}}R_{\Lambda_M}\otimes(R_{\Lambda_{M^\rmt}}+R_{\Lambda_{M^\rmt\cup\Psi'}})
=\sqrt 2\sum_{M\subset Z_\rmI,\ |M|\text{ even}}R_{\Lambda_M}\otimes R^{(1)}_{\Lambda_{M^\rmt}}
\end{align*}
where $\Psi'=\binom{2m+1}{-}$.
Let $\widetilde f\colon\calv_Z\rightarrow\calv_Z$ be the identity map.
Note that now $Z=Z_{(m)}$ and $Z'^{(1)}=Z'_{(m)}$ and
by Lemma~\ref{0302}, we have
\begin{align*}
\widetilde f\otimes\widetilde f'\left(\sum_{M\subset Z_\rmI,\ |M|\text{ even}}
\rho_{\Lambda_M}\otimes\rho^{(1)}_{\Lambda_{M^\rmt}}\right)
&= \sum_{(\Lambda,\Lambda')\in\calb_{Z,Z'^{(1)}}}\rho_\Lambda\otimes\rho_{\Lambda'} \\
&= \sum_{(\Lambda,\Lambda')\in\calb_{Z,Z'^{(1)}}} R_\Lambda\otimes R_{\Lambda'} \\
&= \widetilde f\otimes\widetilde f'\left(\sum_{M\subset Z_\rmI,\ |M|\text{ even}}
R_{\Lambda_M}\otimes R^{(1)}_{\Lambda_{M^\rmt}}\right).
\end{align*}
Because $\widetilde f\otimes\widetilde f'$ is an isometry, the lemma is proved.
\end{proof}

\subsection{The regular cases}
Let $(\bfG,\bfG')=(\Sp_{2n},\rmO^\pm_{2n'})$, and let
\begin{equation}\label{0311}
Z=\binom{a_1,a_2,\ldots,a_{m+1}}{b_1,b_2,\ldots,b_m},\qquad
Z'=\binom{c_1,c_2,\ldots,c_{m'}}{d_1,d_2,\ldots,d_{m'}}
\end{equation}
be special symbols of rank $n,n'$ and of defect $1,0$ respectively.
For any fixed symbol $\Lambda_1\in\cals_Z$ and any fixed symbol $\Lambda'_1\in\cals^\epsilon_{Z'}$,
we define
\begin{align}\label{0324}
\begin{split}
D^\epsilon_{\Lambda_1} &=\bigl\{\,\Lambda'\in\cals_{Z'}^\epsilon
\mid (\Lambda_1,\Lambda')\in\calb^\epsilon_{Z,Z'}\,\bigr\}, \\
D_{\Lambda'_1} &=\bigl\{\,\Lambda\in\cals_Z\mid (\Lambda,\Lambda'_1)\in\calb^\epsilon_{Z,Z'}\,\bigr\}
\end{split}
\end{align}
where $\epsilon=+$ or $-$.
For any fixed symbol $\Lambda'_1\in D^\epsilon_\Lambda$, it is known that (\cf.~\cite{pan-uniform} proposition~6.9)
\[
D^\epsilon_\Lambda=\Lambda'_1+D^+_Z=\{\,\Lambda'_1+\Lambda'\mid\Lambda'\in D^+_Z\,\}.
\]
Similarly, for any fixed symbol $\Lambda_1\in D_{\Lambda'}$, we have
\[
D_{\Lambda'}=\Lambda_1+D_{Z'}=\{\,\Lambda_1+\Lambda\mid\Lambda\in D_{Z'}\,\}.
\]

A special symbol $Z$ is called \emph{regular} if $Z=Z_\rmI$, i.e., $Z^*\cap Z_*=\emptyset$.
The relation $\calb_{Z,Z'}$ for special symbols $Z,Z'$ is called \emph{regular} if
\begin{itemize}
\item both $Z,Z'$ are regular;

\item $D^+_Z=\{Z'\}$ and $D_{Z'}=\{Z\}$.
\end{itemize}

\begin{lem}\label{0303}
Let $Z,Z'$ be special symbols of rank $n,n'$ and of defect $1,0$ respectively.
Suppose that $\calb_{Z,Z'}$ is regular.
Then
\[
\sum_{(\Lambda,\Lambda')\in\calb_{Z,Z'}}\rho_\Lambda\otimes\rho_{\Lambda'}
=\sum_{(\Lambda,\Lambda')\in\calb_{Z,Z'}}R_\Lambda\otimes R_{\Lambda'}.
\]
\end{lem}
\begin{proof}
Write $Z,Z'$ as in (\ref{0311}) for some integers $m,m'$, and suppose that $\calb_{Z,Z'}$ is regular.
Because now $\calb_{Z,Z'}\neq\emptyset$,
we know that either $m'=m$ or $m'=m+1$ (\cf.~\cite{pan-uniform} lemma~2.20).
Let $(\bfG,\bfG')=(\Sp_{2n},\rmO^\pm_{2n'})$ and $(\bfG^{(1)},\bfG'^{(1)})=(\Sp_{2m(m+1)},\rmO^\pm_{2m'^2})$.
\begin{enumerate}
\item First suppose that $m'=m$.
Define
\begin{align*}
\theta\colon\{c_1,c_2,\ldots,c_m\}\cup\{d_1,d_2,\ldots,d_m\}
&\rightarrow\{a_1,a_2,\ldots,a_{m+1}\}\cup\{b_1,b_2,\ldots,b_m\} \\
c_i &\mapsto b_i \\
d_i &\mapsto a_{i+1}.
\end{align*}
Note that $a_1$ is not in the image of $\theta$.
For $N\subset Z'_\rmI$, we see that both $\theta(N),\binom{a_1}{-}\cup\theta(N)$ are sub-symbols
of $Z_\rmI$
where $\theta(N)=\{\,\theta(x)\mid x\in N\,\}$.
Then $\theta$ induces a mapping $\bar\theta\colon\cals_{Z'}\rightarrow\cals_Z$ given by
\[
\bar\theta(\Lambda_{N})=\begin{cases}
\Lambda_{\theta(N)}, & \text{if $N\subset Z'_\rmI$, $|N|$ even};\\
\Lambda_{\binom{a_1}{-}\cup\theta(N)}, & \text{if $N\subset Z'_\rmI$, $|N|$ odd}.
\end{cases}
\]
From \cite{pan-uniform} lemma~6.20 and lemma 8.5, we know that
\begin{equation}\label{0313}
\calb_{Z,Z'}=\left\{\,(\bar\theta(\Lambda'),\Lambda')\mid\Lambda'\in\cals_{Z'}\,\right\}.
\end{equation}

\item Suppose that $m'=m+1$.
Define
\begin{align*}
\theta\colon\{a_1,a_2,\ldots,a_{m+1}\}\cup\{b_1,b_2,\ldots,b_m\}
&\rightarrow\{c_1,c_2,\ldots,c_{m+1}\}\cup\{d_1,d_2,\ldots,d_{m+1}\} \\
a_i &\mapsto d_i \\
b_i &\mapsto c_{i+1}.
\end{align*}
Note that $c_1$ is not in the image of $\theta$.
For $M\subset Z_\rmI$ such that $|M|$ is even,
we see that $\theta(M)$ and $\theta(M)\cup\Psi'$ are sub-symbols of $Z'$
where $\Psi'=\binom{c_1}{-}$.
The $\theta$ induces two mappings $\bar\theta^+,\bar\theta^-\colon\cals_Z\rightarrow\cals_{Z'}$ by
\begin{align*}
\bar\theta^+(\Lambda_M) &=\Lambda_{\theta(M)}\in\cals^+_{Z'}; \\
\bar\theta^-(\Lambda_M) &=\Lambda_{\theta(M)\cup\Psi'}\in\cals^-_{Z'}
\end{align*}
By \cite{pan-uniform} lemma~6.20 and lemma~8.5, we have
\begin{equation}\label{0314}
\calb_{Z,Z'}=\left\{\,(\Lambda,\bar\theta^+(\Lambda)),(\Lambda,\bar\theta^-(\Lambda))
\mid\Lambda\in\cals_Z\,\right\}.
\end{equation}
\end{enumerate}

Now define a bijection
\begin{align*}
f\colon\{a_1,a_2,\ldots,a_{m+1},b_1,b_2,\ldots,b_m\} &\rightarrow\{2m,2m-1,\ldots,0\} \\
a_i &\mapsto 2m+2-2i; \\
b_i &\mapsto 2m+1-2i.
\end{align*}
Then $f$ induces a bijection
$\bar f\colon\cals_Z\rightarrow\cals_{Z_{(m)}}$ given by $\Lambda_M\mapsto\Lambda_{f(M)}$ for $M\subset Z_\rmI$.
Because $|M_1\cap M_2|=|f(M_1)\cap f(M_2)|$, we have
$\langle\rho_{\Lambda_{M_1}},R_{\Lambda_{M_2}}\rangle_\bfG
=\langle\rho_{\bar f(\Lambda_{M_1})},R_{\bar f(\Lambda_{M_2})}\rangle_{\bfG^{(1)}}$.
Similarly, we define
\begin{align*}
f'\colon\{c_1,c_2,\ldots,c_{m'},d_1,d_2,\ldots,d_{m'}\} &\rightarrow\{2m'-1,2m'-2,\ldots,0\} \\
c_i &\mapsto 2m'+1-2i; \\
d_i &\mapsto 2m'-2i.
\end{align*}
Then we have  a bijection
$\bar f'\colon\cals_{Z'}\rightarrow\cals_{Z'_{(m)}}$ given by $\Lambda_{M'}\mapsto\Lambda_{f'(M')}$.
Note that
\[
|M_1'||M_2'|+|M'_1\cap M'_2|=|f'(M_1')||f'(M_2')|+|f'(M'_1)\cap f'(M'_2)|,
\]
so we have
$\langle\rho_{\Lambda'_1},R_{\Lambda'_2}\rangle_{\bfG'}
=\langle\rho_{\bar f'(\Lambda'_1)},R_{\bar f'(\Lambda'_2)}\rangle_{\bfG'^{(1)}}$.
By (\ref{0313}) and (\ref{0314}), we see that $(\Lambda,\Lambda')\in\calb_{Z,Z'}$ if and only if
$(\bar f(\Lambda),\bar f'(\Lambda'))\in\calb_{Z_{(m)},Z'_{(m')}}$
for $\Lambda\in\cals_Z$ and $\Lambda'\in\cals_{Z'}$.

Now $\bar f\times\bar f'$ induces an isometry
\[
\widetilde f\otimes\widetilde f'\colon\calv_Z\otimes\calv_{Z'}\rightarrow\calv_{Z_{(m)}}\otimes\calv_{Z'_{(m')}}
\]
such that $\widetilde f(\rho_\Lambda)=\rho_{\bar f(\Lambda)}$,
$\widetilde f(R_\Lambda)=R_{\bar f(\Lambda)}$,
$\widetilde f'(\rho_{\Lambda'})=\rho_{\bar f'(\Lambda')}$,
and $\widetilde f'(R_{\Lambda'})=R_{\bar f'(\Lambda')}$.
Therefore, by Lemma~\ref{0302} and Lemma~\ref{0307}, we have
\begin{align*}
\widetilde f\otimes\widetilde f'\left(\sum_{(\Lambda,\Lambda')\in\calb_{Z,Z'}}\rho_\Lambda\otimes\rho_{\Lambda'}\right)
&=\sum_{(\bar f(\Lambda),\bar f'(\Lambda'))\in\calb_{Z_{(m)},Z'_{(m')}}}
\rho_{\bar f(\Lambda)}\otimes\rho_{\bar f'(\Lambda')} \\
&=\sum_{(\bar f(\Lambda),\bar f'(\Lambda'))\in\calb_{Z_{(m)},Z'_{(m')}}}
R_{\bar f(\Lambda)}\otimes R_{\bar f'(\Lambda')} \\
&=\widetilde f\otimes\widetilde f'\left(\sum_{(\Lambda,\Lambda')\in\calb_{Z,Z'}}R_\Lambda\otimes R_{\Lambda'}\right).
\end{align*}
Therefore the lemma is proved.
\end{proof}

\subsection{General cases}\label{0310}
Let $(\bfG,\bfG')=(\Sp_{2n},\rmO^\pm_{2n'})$, and let
\[
Z=\binom{a_1,a_2,\ldots,a_{m+1}}{b_1,b_2,\ldots,b_m},\qquad
Z'=\binom{c_1,c_2,\ldots,c_{m'}}{d_1,d_2,\ldots,d_{m'}}
\]
be special symbols of rank $n,n'$ and of defect $1,0$ respectively.
We assume that $\calb_{Z,Z'}\neq\emptyset$ and so we have $m'=m+1$ or $m'=m$ (\cf.~\cite{pan-uniform} lemma~2.20).

A pair $\binom{a_k}{b_l}\subset Z$ is called an \emph{identical pair} if $a_k=b_l$.
A pair $\binom{a_k}{b_l}\subset Z_\rmI$ is called an \emph{consecutive pair} if there is no
entry $x$ of $Z$ such that either $a_k<x<b_l$ or $b_l<x<a_k$.
Similar definitions apply to $Z'$.
A consecutive pair $\Psi=\binom{a_k}{b_l}\subset Z_\rmI$ is said to be in the \emph{core}
of $\calb_{Z,Z'}$ if $\Lambda_\Psi\in D_{Z'}$ where $\Lambda_\Psi$ is given as in (\ref{0210});
similarly a consecutive pair $\Psi'=\binom{c_{l'}}{d_{k'}}\subset Z'_\rmI$ is said to be in the
\emph{core} of $\calb_{Z,Z'}$ if $\Lambda_{\Psi'}\in D^+_Z$
where $D_{Z'}$ and $D^+_Z$ are given in (\ref{0324}).

We consider sets of pairs in $Z,Z'$ respectively according the following order
\begin{equation}\label{0814}
\textstyle
\bigl\{\binom{a_{m+1}}{b_m},\binom{c_{m'}}{d_{m'}}\bigr\},
\bigl\{\binom{a_m}{b_m},\binom{c_{m'}}{d_{m'-1}}\bigr\},
\bigl\{\binom{a_m}{b_{m-1}},\binom{c_{m'-1}}{d_{m'-1}}\bigr\},
\bigl\{\binom{a_{m-1}}{b_{m-1}},\binom{c_{m'-1}}{d_{m'-2}}\bigr\},\ldots.
\end{equation}
Suppose that $\left\{\Psi:=\binom{a_k}{b_l},\Psi':=\binom{c_{l'}}{d_{k'}}\right\}$ is the first set
of pairs in the above order such that at least one of the two pairs in the set is either
an identical pair or in the core of $\calb_{Z,Z'}$.
If $m=0,m'=1$, by convention, let $\Psi=\emptyset,\Psi'=\binom{c_1}{d_1}$.
Such a set $\{\Psi,\Psi'\}$ exists if $\calb_{Z,Z'}$ is not regular.
Note that if $m=m'=0$, then $Z=\binom{a_1}{-},Z'=\binom{-}{-}$ and $\calb_{Z,Z'}$ is always regular.
Now we have the following situations:
\begin{enumerate}
\item[(I)] If each of $\binom{a_k}{b_l}$ and $\binom{c_{l'}}{d_{k'}}$ is either an identical pair or in
the core of $\calb_{Z,Z'}$, then we define
\begin{align*}
Z^{(1)} &= \binom{a_1-1,a_2-1,\ldots,a_{k-1}-1,a_{k+1},\ldots,a_{m+1}}
{b_1-1,b_2-1,\ldots,b_{l-1}-1,b_{l+1},\ldots,b_m}, \\
Z'^{(1)} &=\binom{c_1-1,c_2-1,\ldots,c_{l'-1}-1,c_{l'+1},\ldots,c_{m'}}
{d_1-1,d_2-1,\ldots,d_{k'-1}-1,d_{k'+1},\ldots,d_{m'}}.
\end{align*}
Define a bijective mapping $f\colon Z\smallsetminus\binom{a_k}{b_l}\rightarrow Z^{(1)}$ given by
\begin{equation}
a_i\mapsto\begin{cases}
a_i-1, & \text{if $1\leq i<k$};\\
a_i, & \text{if $k<i\leq m+1$},
\end{cases}\qquad
b_j\mapsto\begin{cases}
b_j-1, & \text{if $1\leq j<l$};\\
b_j, & \text{if $l<j\leq m$}.
\end{cases}
\end{equation}
Similarly, we define a bijective mapping $f'\colon Z'\smallsetminus\binom{c_{l'}}{d_{k'}}\rightarrow Z'^{(1)}$
given by
\begin{equation}
c_i\mapsto\begin{cases}
c_i-1, & \text{if $1\leq i<l'$};\\
c_i, & \text{if $l'<i\leq m'$},
\end{cases}\qquad
d_j\mapsto\begin{cases}
d_j-1, & \text{if $1\leq j<k'$};\\
d_j, & \text{if $k'<j\leq m'$}.
\end{cases}
\end{equation}

\item[(II)] If $\binom{a_k}{b_l}$ is either an identical pair or in the core, but $\binom{c_{l'}}{d_{k'}}$ is neither,
then we define
\begin{align*}
Z^{(1)} &= \binom{a_1-1,a_2-1,\ldots,a_{k-1}-1,d_{k'+1},\ldots,d_{m'}}
{b_1-1,b_2-1,\ldots,b_{l-1}-1,c_{l'+1},\ldots,c_m'}, \\
Z'^{(1)} &=\binom{c_1,c_2,\ldots,c_{l'},b_{l+1},\ldots,b_{m}}
{d_1,d_2,\ldots,d_{k'},a_{k+1},\ldots,a_{m+1}}.
\end{align*}
Note that this case happens only when $m'=m$ (\cf.~\cite{pan-uniform} lemma~7.2).
Define a bijective mapping $f\colon Z\smallsetminus\binom{a_k}{b_l}\rightarrow Z^{(1)}$ given by
\begin{equation}
a_i\mapsto\begin{cases}
a_i-1, & \text{if $1\leq i<k$};\\
d_{i-1}, & \text{if $k<i\leq m+1$},
\end{cases}\qquad
b_j\mapsto\begin{cases}
b_j-1, & \text{if $1\leq j<l$};\\
c_j, & \text{if $l<j\leq m$}.
\end{cases}
\end{equation}
Similarly, define a bijective mapping $f'\colon Z'\rightarrow Z'^{(1)}$ given by
\begin{equation}
c_i\mapsto\begin{cases}
c_i, & \text{if $1\leq i\leq l'$};\\
b_i, & \text{if $l'<i\leq m'$},
\end{cases}\qquad
d_j\mapsto\begin{cases}
d_j, & \text{if $1\leq j\leq k'$};\\
a_{j+1}, & \text{if $k'<j\leq m'$}.
\end{cases}
\end{equation}

\item[(III)] If $\binom{c_{l'}}{d_{k'}}$ is either an identical pair or in the core, but $\binom{a_k}{b_l}$
is neither, then we define
\begin{align*}
Z^{(1)} &= \binom{a_1,a_2,\ldots,a_k,d_{k'+1},\ldots,d_{m'}}
{b_1,b_2,\ldots,b_l,c_{l'+1},\ldots,c_{m'}}, \\
Z'^{(1)} &=\binom{c_1-1,c_2-1,\ldots,c_{l'-1}-1,b_{l+1},\ldots,b_m}
{d_1-1,d_2-1,\ldots,d_{k'-1}-1,a_{k+1},\ldots,a_{m+1}}.
\end{align*}
Note that this case happens only when $m'=m+1$ (\cf.~\cite{pan-uniform} lemma~7.2).
Define a bijective mapping $f\colon Z\rightarrow Z^{(1)}$ given by
\begin{equation}
a_i\mapsto\begin{cases}
a_i, & \text{if $1\leq i\leq k$};\\
d_i, & \text{if $k<i\leq m+1$},
\end{cases}\qquad
b_j\mapsto\begin{cases}
b_j, & \text{if $1\leq j\leq l$};\\
c_{j+1}, & \text{if $l<j\leq m$}.
\end{cases}
\end{equation}
Similarly, define a bijective mapping $f'\colon Z'\smallsetminus\binom{c_{l'}}{d_{k'}}\rightarrow Z'^{(1)}$ given by
\begin{equation}
c_i\mapsto\begin{cases}
c_i-1, & \text{if $1\leq i<l'$};\\
b_{i-1}, & \text{if $l'<i\leq m'$},
\end{cases}\qquad
d_j\mapsto\begin{cases}
d_j-1, & \text{if $1\leq j<k'$};\\
a_j, & \text{if $k'<j\leq m'$}.
\end{cases}
\end{equation}
\end{enumerate}
For all three cases (I),(II),(III),
let $\bfG^{(1)}$ be the symplectic group $\Sp_{2n^{(1)}}$ such that $Z^{(1)}\in\cals_{\bfG^{(1)}}$,
and let $\bfG'^{(1)}$ be the orthogonal group $\rmO^\pm_{2n'^{(1)}}$ such that
$Z'^{(1)}\in\cals_{\bfG'^{(1)}}$.
Moreover, we also define
\begin{align*}
\cals^{(1)}_Z &= \begin{cases}
\{\,\Lambda_M\mid M\subset Z_\rmI\smallsetminus\binom{a_k}{b_l},\ |M|\ \text{even}\,\},
& \text{if $\binom{a_k}{b_l}$ is an identical pair or in the core};\\
\{\,\Lambda_M\mid M\subset Z_\rmI,\ |M|\ \text{even}\,\},
& \text{otherwise};
\end{cases} \\
\cals^{(1)}_{Z'} &= \begin{cases}
\{\,\Lambda_{M'}\mid M'\subset Z'_\rmI\smallsetminus\binom{c_{l'}}{d_{k'}}\,\},
& \text{if $\binom{c_{l'}}{d_{k'}}$ is an identical pair or in the core}; \\
\{\,\Lambda_{M'}\mid M'\subset Z'_\rmI\,\}
& \text{otherwise}
\end{cases} \\
\calb^{(1)}_{Z,Z'} &=\calb_{Z,Z'}\cap(\cals^{(1)}_Z\times\cals^{(1)}_{Z'}).
\end{align*}
Note that ${\rm size}(Z)=(m+1,m)$ and ${\rm size}(Z')=(m',m')$.
From the above construction, it is not difficult to see that
$Z^{(1)},Z'^{(1)}$ are special symbols and
\[
{\rm size}(Z^{(1)})=\begin{cases}
(m,m-1) \\
(m,m-1) \\
(m+1,m)
\end{cases}\text{ and}\quad
{\rm size}(Z'^{(1)})=\begin{cases}
(m'-1,m'-1), & \text{for Case (I)};\\
(m',m'), & \text{for Case (II)};\\
(m'-1,m'-1), & \text{for Case (III)}.
\end{cases}
\]
In particular,
${\rm def}(Z^{(1)})={\rm def}(Z)=1$, ${\rm def}(Z'^{(1)})={\rm def}(Z')=0$.
Moreover,
\begin{align*}
\deg(Z^{(1)}) &=\begin{cases}
\deg(Z), & \text{if $\binom{a_k}{b_l}$ is not in the core};\\
\deg(Z)-1, & \text{if $\binom{a_k}{b_l}$ is in the core},
\end{cases} \\
\deg(Z'^{(1)}) &=\begin{cases}
\deg(Z'), & \text{if $\binom{c_{l'}}{d_{k'}}$ is not in the core};\\
\deg(Z')-1, & \text{if $\binom{c_{l'}}{d_{k'}}$ is in the core}.
\end{cases}
\end{align*}

Now we define $\bar f,\bar f'$ according the following situations:
\begin{enumerate}
\item Suppose that $\binom{a_k}{b_l}$ is not in the core.
Then $f$ induces a bijection $\bar f\colon \cals_Z^{(1)}\rightarrow\cals_{Z^{(1)}}$
by $\bar f(\Lambda_M)=\Lambda_{f(M)}$ for $M\subset Z_\rmI$ and $|M|$ even.
For this case, let $\calv_Z^{(1)}=\calv_Z$, $\rho_\Lambda^{(1)}=\rho_\Lambda$ and $R_\Lambda^{(1)}=R_\Lambda$
for $\Lambda\in\cals_Z$.

\item Suppose that $\Psi:=\binom{a_k}{b_l}$ is in the core.
Then $f$ induces a bijection $\bar f\colon\cals_Z^{(1)}\rightarrow\cals_{Z^{(1)}}$ by $\bar f(\Lambda_M)=\Lambda_{f(M)}$
for $M\subset Z_\rmI\smallsetminus\Psi$ and $|M|$ even.
Define
\begin{equation}\label{0320}
\rho^{(1)}_{\Lambda_M}:=\frac{1}{\sqrt 2}(\rho_{\Lambda_M}+\rho_{\Lambda_{M\cup\Psi}}),\qquad
R^{(1)}_{\Lambda_M}:=\frac{1}{\sqrt 2}(R_{\Lambda_M}+R_{\Lambda_{M\cup\Psi}}),
\end{equation}
and let $\calv^{(1)}_Z$ be the subspace spanned by $\rho^{(1)}_{\Lambda_M}$ for
$M\subset Z_\rmI\smallsetminus\Psi$ and $|M|$ even.

\item Suppose that $\binom{c_{l'}}{d_{k'}}$ is not in the core.
Then $f'$ induces a bijection $\bar f'\colon \cals_{Z'}^{(1)}\rightarrow\cals_{Z'^{(1)}}$
by $\bar f'(\Lambda_{M'})=\Lambda_{f'(M')}$ for $M'\subset Z'_\rmI$.
For this case, let $\calv^{(1)}_{Z'}=\calv_{Z'}$, $\rho^{(1)}_{\Lambda'}=\rho_{\Lambda'}$ and
$R^{(1)}_{\Lambda'}=R_{\Lambda'}$ for $\Lambda'\in\cals_{Z'}$.

\item Suppose that $\Psi':=\binom{c_{l'}}{d_{k'}}$ is in the core.
Then $f'$ induces a bijection $\bar f'\colon\cals_{Z'}^{(1)}\rightarrow\cals_{Z'^{(1)}}$
by $\bar f'(\Lambda_M')=\Lambda_{f'(M')}$
for $M'\subset Z'_\rmI\smallsetminus\Psi'$.
Define
\begin{equation}\label{0321}
\rho^{(1)}_{\Lambda_{M'}}:=\frac{1}{\sqrt 2}(\rho_{\Lambda_{M'}}+\rho_{\Lambda_{M'\cup\Psi'}}),\qquad
R^{(1)}_{\Lambda_{M'}}:=\frac{1}{\sqrt 2}(R_{\Lambda_{M'}}+R_{\Lambda_{M'\cup\Psi'}}),
\end{equation}
and let $\calv^{(1)}_{Z'}$ be the subspace spanned by $\rho^{(1)}_{\Lambda_{M'}}$ for
$M'\subset Z'_\rmI\smallsetminus\Psi'$.
\end{enumerate}

\begin{lem}\label{0317}
Keep the above notations.
Then
\begin{align*}
\sum_{(\Lambda,\Lambda')\in\calb_{Z,Z'}}\rho_\Lambda\otimes\rho_{\Lambda'}
&= C \sum_{(\Lambda,\Lambda')\in\calb^{(1)}_{Z,Z'}}\rho^{(1)}_\Lambda\otimes\rho^{(1)}_{\Lambda'} \\
\sum_{(\Lambda,\Lambda')\in\calb_{Z,Z'}}R_\Lambda\otimes R_{\Lambda'}
&= C \sum_{(\Lambda,\Lambda')\in\calb^{(1)}_{Z,Z'}}R^{(1)}_\Lambda\otimes R^{(1)}_{\Lambda'}
\end{align*}
where
\[
C=\begin{cases}
1, & \text{if none of the pairs $\binom{a_k}{b_l},\binom{c_{l'}}{d_{k'}}$ is in the core};\\
\sqrt 2, & \text{if exactly one of the pairs $\binom{a_k}{b_l},\binom{c_{l'}}{d_{k'}}$ is in the core};\\
2, & \text{if both of the pairs $\binom{a_k}{b_l},\binom{c_{l'}}{d_{k'}}$ are in the core}.
\end{cases}
\]
\end{lem}
\begin{proof}
From \cite{pan-uniform} lemma~7.10, lemma~7.15, and lemma~8.8, we know that
\begin{align*}
\sum_{(\Lambda,\Lambda')\in\calb^\epsilon_{Z,Z'}}\rho_\Lambda\otimes\rho_{\Lambda'}
&= C \sum_{(\Lambda,\Lambda')\in\calb^{\epsilon,(1)}_{Z,Z'}}\rho^{(1)}_\Lambda\otimes\rho^{(1)}_{\Lambda'} \\
\sum_{(\Lambda,\Lambda')\in\calb^\epsilon_{Z,Z'}}R_\Lambda\otimes R_{\Lambda'}
&= C \sum_{(\Lambda,\Lambda')\in\calb^{\epsilon,(1)}_{Z,Z'}}R^{(1)}_\Lambda\otimes R^{(1)}_{\Lambda'}
\end{align*}
where $\calb^{\epsilon,(1)}_{Z,Z'}=\calb^{(1)}_{Z,Z'}\cap\calb^\epsilon_{Z,Z'}$ for $\epsilon=+$ or $-$.
Moreover, we know that
\[
\calb_{Z,Z'}=\calb^+_{Z,Z'}\cup\calb^-_{Z,Z'}\quad\text{and}\quad
\calb^{(1)}_{Z,Z'}=\calb^{+,(1)}_{Z,Z'}\cup\calb^{-,(1)}_{Z,Z'}.
\]
Therefore, the lemma is proved.
\end{proof}

\begin{lem}\label{0318}
Keep the above notations.
Then $(\Lambda,\Lambda')\in\calb^{(1)}_{Z,Z'}$ if and only if\/
$(\bar f(\Lambda),\bar f'(\Lambda'))\in\calb_{Z^{(1)},Z'^{(1)}}$.
\end{lem}
\begin{proof}
From \cite{pan-uniform} lemma~7.9, lemma~7.14 and lemma~8.7, we know that
$(\Lambda,\Lambda')\in\calb^{\epsilon,(1)}_{Z,Z'}$ if and only if
$(\bar f(\Lambda),\bar f'(\Lambda'))\in\calb^\epsilon_{Z^{(1)},Z'^{(1)}}$ for $\epsilon=+$ or $-$.
Now we have $\calb^{(1)}_{Z,Z'}=\calb^{+,(1)}_{Z,Z'}\cup\calb^{-,(1)}_{Z,Z'}$.
Hence the lemma follows immediately.
\end{proof}

\begin{lem}\label{0316}
Keep the above notations.
Then
\[
\langle \rho^{(1)}_{\Lambda_1},R^{(1)}_{\Lambda_2}\rangle_\bfG
=\langle \rho_{\bar f(\Lambda_1)},R_{\bar f(\Lambda_2)}\rangle_{\bfG^{(1)}}\quad\text{and}\quad
\langle \rho^{(1)}_{\Lambda'_1},R^{(1)}_{\Lambda'_2}\rangle_{\bfG'}
=\langle \rho_{\bar f'(\Lambda'_1)},R_{\bar f'(\Lambda'_2)}\rangle_{\bfG'^{(1)}}
\]
for any $\Lambda_1,\Lambda_2\in\cals^{(1)}_Z$ and any $\Lambda'_1,\Lambda'_2\in\cals^{(1)}_{Z'}$.
\end{lem}
\begin{proof}
Write $\Lambda_1=\Lambda_{M_1}$, $\Lambda_2=\Lambda_{M_2}$, $\Lambda'_1=\Lambda_{M'_1}$ and
$\Lambda'_2=\Lambda_{M'_2}$ for some $M_1,M_2,M'_1, M'_2$.
Now we have the following situations.
\begin{enumerate}
\item Suppose that $\binom{a_k}{b_l}$ is not in the core of $\calb_{Z,Z'}$.
Then $\rho^{(1)}_{\Lambda_{M_1}}=\rho_{\Lambda_{M_1}}$,
$R^{(1)}_{\Lambda_{M_2}}=R_{\Lambda_{M_2}}$,
$|M_1\cap M_2|=|f(M_1)\cap f(M_2)|$,
$\deg(Z^{(1)})=\deg(Z)$.
From (\ref{0203}) and (\ref{0205}), we have $\langle \rho^{(1)}_{\Lambda_1},R^{(1)}_{\Lambda_2}\rangle_\bfG
=\langle \rho_{\bar f(\Lambda_1)},R_{\bar f(\Lambda_2)}\rangle_{\bfG^{(1)}}$ immediately.

\item Suppose that $\Psi:=\binom{a_k}{b_l}$ is in the core of $\calb_{Z,Z'}$.
Now $M_1,M_2\subset Z_\rmI\smallsetminus\Psi$,
\begin{align*}
|M_1\cap M_2|=|M_1\cap(M_2\cup\Psi)|=|(M_1\cup\Psi)\cap M_2| &=|f(M_1)\cap f(M_2)|, \\
|(M_1\cup\Psi)\cap(M_2\cup\Psi)| &=|f(M_1)\cap f(M_2)|+2,
\end{align*}
and $\deg(Z^{(1)})=\deg(Z)-1$.
Therefore,
\begin{align*}
&\langle\rho^{(1)}_{\Lambda_{M_1}},R^{(1)}_{\Lambda_{M_2}}\rangle_\bfG \\
&= \frac{1}{2}\left[\langle\rho_{\Lambda_{M_1}},R_{\Lambda_{M_2}}\rangle_\bfG
+\langle\rho_{\Lambda_{M_1\cup\Psi}},R_{\Lambda_{M_2}}\rangle_\bfG
+\langle\rho_{\Lambda_{M_1}},R_{\Lambda_{M_2\cup\Psi}}\rangle_\bfG
+\langle\rho_{\Lambda_{M_1\cup\Psi}},R_{\Lambda_{M_2\cup\Psi}}\rangle_\bfG\right] \\
&= 2\,\langle\rho_{\Lambda_{M_1}},R_{\Lambda_{M_2}}\rangle_\bfG
=\langle\rho_{\bar f(\Lambda_{M_1})},R_{\bar f(\Lambda_{M_2})}\rangle_{\bfG^{(1)}}.
\end{align*}

\item Suppose that $\binom{c_{l'}}{d_{k'}}$ is not in the core of $\calb_{Z,Z'}$.
Now $\rho^{(1)}_{\Lambda_{M'_1}}=\rho_{\Lambda_{M'_1}}$,
$R^{(1)}_{\Lambda_{M'_2}}=R_{\Lambda_{M'_2}}$,
$|M'_1||M'_2|=|f(M'_1)||f(M'_2)|$,
$|M'_1\cap M'_2|=|f(M'_1)\cap f(M'_2)|$, and
$\deg(Z'^{(1)})=\deg(Z')$.
From (\ref{0204}) and (\ref{0206}),
we have $\langle\rho^{(1)}_{\Lambda'_1},R^{(1)}_{\Lambda'_2}\rangle_{\bfG'}
=\langle\rho_{\bar f'(\Lambda'_1)},R_{\bar f'(\Lambda'_2)}\rangle_{\bfG'^{(1)}}$ immediately.

\item Suppose that $\Psi':=\binom{c_{l'}}{d_{k'}}$ is in the core of $\calb_{Z,Z'}$.
Now $M'_1,M'_2\subset Z'_\rmI\smallsetminus\Psi'$.
Then
\begin{align*}
|M'_1\cap M'_2|=|M'_1\cap(M'_2\cup\Psi')|=|(M'_1\cup\Psi')\cap M'_2| &=|f'(M'_1)\cap f'(M'_2)|, \\
|(M'_1\cup\Psi')\cap(M'_2\cup\Psi')| &=|f'(M'_1)\cap f'(M'_2)|+2.
\end{align*}
Moreover, we have
\begin{align*}
|M'_1||M'_2|\equiv|M'_1||M'_2\cup\Psi'|\equiv|M'_1\cup\Psi'||M'_2|
&\equiv |M'_1\cup\Psi'||M'_2\cup\Psi'| \\
&\equiv|f'(M'_1)||f'(M'_2)| \pmod 2
\end{align*}
and $\deg(Z'^{(1)})=\deg(Z')-1$.
Therefore, as in (1), we have
\[
\langle\rho^{(1)}_{\Lambda_{M'_1}},R^{(1)}_{\Lambda_{M'_2}}\rangle_{\bfG'}
= 2\,\langle\rho_{\Lambda_{M'_1}},R_{\Lambda_{M'_2}}\rangle_{\bfG'}
= \langle\rho_{\bar f'(\Lambda_{M'_1})},R_{\bar f'(\Lambda_{M'_2})}\rangle_{\bfG'^{(1)}}.
\]
\end{enumerate}
\end{proof}

Now the bijection $\bar f\times\bar f'\colon\cals^{(1)}_Z\times\cals^{(1)}_{Z'}\rightarrow\cals_{Z^{(1)}}\times\cals_{Z'^{(1)}}$
induces an isometry
\[
\widetilde f\otimes\widetilde f'\colon\calv_Z^{(1)}\otimes\calv^{(1)}_{Z'}
\rightarrow\calv_{Z^{(1)}}\otimes\calv_{Z'^{(1)}}\quad\text{by}\quad
\rho^{(1)}_{\Lambda}\otimes\rho^{(1)}_{\Lambda'}\mapsto\rho_{\bar f(\Lambda)}\otimes\rho_{\bar f'(\Lambda')}.
\]

\begin{lem}\label{0319}
Keep the above notations.
Then $\widetilde f(R_\Lambda^{(1)})=R_{\bar f(\Lambda)}$ and
$\widetilde f'(R_{\Lambda'}^{(1)})=R_{\bar f'(\Lambda')}$
for any $\Lambda\in\cals_Z^{(1)}$ and any $\Lambda'\in\cals_{Z'}^{(1)}$.
\end{lem}
\begin{proof}
Now $\bigl\{\,\rho^{(1)}_{\Lambda}\otimes\rho^{(1)}_{\Lambda'}\mid\Lambda\in\cals_Z^{(1)},\ \Lambda'\in\cals^{(1)}_{Z'}\,\bigr\}$ is an orthonormal basis for $\calv^{(1)}_Z\otimes\calv^{(1)}_{Z'}$, and
$\bigl\{\,\rho_{\Lambda^{(1)}}\otimes\rho_{\Lambda'^{(1)}}\mid\Lambda^{(1)}\in\cals_{Z^{(1)}},\ \Lambda'^{(1)}\in\cals^{(1)}_{Z'}\,\bigr\}$
is an orthonormal basis for $\calv_{Z^{(1)}}\otimes\calv_{Z'^{(1)}}$.
Moreover,
$\widetilde f\otimes\widetilde f'\colon \calv^{(1)}_Z\otimes\calv^{(1)}_{Z'}\rightarrow\calv_{Z^{(1)}}\otimes\calv_{Z'^{(1)}}$
given by $\rho^{(1)}_{\Lambda}\otimes\rho^{(1)}_{\Lambda'}\mapsto\rho_{\bar f(\Lambda)}\otimes\rho_{\bar f'(\Lambda')}$
is an isometry of inner product spaces.
Therefore the lemma follows from Lemma~\ref{0316} immediately.
\end{proof}

\begin{lem}\label{0315}
Keep the above notations.
Then
\begin{align*}
\widetilde f\otimes\widetilde f'\Biggl(\sum_{(\Lambda,\Lambda')\in\calb_{Z,Z'}}\rho_\Lambda\otimes\rho_{\Lambda'}\Biggr)
&=C\sum_{(\Lambda^{(1)},\Lambda'^{(1)})\in\calb_{Z^{(1)},Z'^{(1)}}}\rho_{\Lambda^{(1)}}\otimes\rho_{\Lambda'^{(1)}}; \\
\widetilde f\otimes\widetilde f'\Biggl(\sum_{(\Lambda,\Lambda')\in\calb_{Z,Z'}}R_\Lambda\otimes R_{\Lambda'}\Biggr)
&=C\sum_{(\Lambda^{(1)},\Lambda'^{(1)})\in\calb_{Z^{(1)},Z'^{(1)}}}R_{\Lambda^{(1)}}\otimes R_{\Lambda'^{(1)}}
\end{align*}
where constant $C$ is given as in Lemma~\ref{0317}.
\end{lem}
\begin{proof}
The lemma follows from Lemma~\ref{0317}, Lemma~\ref{0318} and Lemma~\ref{0319} immediately.
\end{proof}

\begin{proof}[Proof of Proposition~\ref{0301}]
Let $Z,Z'$ be special symbols of ranks $n,n'$ and defects $1,0$ respectively.
If $\calb_{Z,Z'}$ is regular, then the proposition is proved in Lemma~\ref{0303}.
Now suppose that $\calb_{Z,Z'}$ is not regular.
From above construction in this subsection, we have
\begin{itemize}
\item special symbols $Z^{(1)},Z'^{(1)}$ of defects $1,0$ respectively,

\item subspaces $\calv^{(1)}_Z\subset\calv_Z$ and $\calv^{(1)}_{Z'}\subset\calv_{Z'}$,

\item an isometry $\widetilde f\otimes\widetilde f'\colon\calv_Z^{(1)}\otimes\calv^{(1)}_{Z'}
\rightarrow\calv_{Z^{(1)}}\otimes\calv_{Z'^{(1)}}$ of inner product spaces.
\end{itemize}
Moreover, by Lemma~\ref{0315}, we have
\begin{align*}
\widetilde f\otimes\widetilde f'\Biggl(\sum_{(\Lambda,\Lambda')\in\calb_{Z,Z'}}\rho_\Lambda\otimes\rho_{\Lambda'}\Biggr)
&=C\sum_{(\Lambda^{(1)},\Lambda'^{(1)})\in\calb_{Z^{(1)},Z'^{(1)}}}\rho_{\Lambda^{(1)}}\otimes\rho_{\Lambda'^{(1)}}; \\
\widetilde f\otimes\widetilde f'\Biggl(\sum_{(\Lambda,\Lambda')\in\calb_{Z,Z'}}R_\Lambda\otimes R_{\Lambda'}\Biggr)
&=C\sum_{(\Lambda^{(1)},\Lambda'^{(1)})\in\calb_{Z^{(1)},Z'^{(1)}}}R_{\Lambda^{(1)}}\otimes R_{\Lambda'^{(1)}}
\end{align*}
for some constant $C$.
Note that ${\rm size}(Z^{(1)})<{\rm size}(Z)$ or ${\rm size}(Z'^{(1)})<{\rm size}(Z')$.
We can repeat the same construction several times until we obtain a regular relation $\calb_{Z^{(t)},Z'^{(t)}}$
and an isometry $\widetilde f_t\otimes\widetilde f'_t\colon\calv_Z^{(t)}\otimes\calv^{(t)}_{Z'}
\rightarrow\calv_{Z^{(t)}}\otimes\calv_{Z'^{(t)}}$.
Note that both two vectors $\sum_{(\Lambda,\Lambda')\in\calb_{Z,Z'}}\rho_\Lambda\otimes\rho_{\Lambda'}$
and $\sum_{(\Lambda,\Lambda')\in\calb_{Z,Z'}}R_\Lambda\otimes R_{\Lambda'}$
are in the subspace $\calv_Z^{(t)}\otimes\calv^{(t)}_{Z'}$.
Then by Lemma~\ref{0315} and Lemma~\ref{0303}, we conclude that
\begin{align*}
\widetilde f_t\otimes\widetilde f'_t\Biggl(\sum_{(\Lambda,\Lambda')\in\calb_{Z,Z'}}\rho_\Lambda\otimes\rho_{\Lambda'}\Biggr)
&=C_t\sum_{(\Lambda^{(t)},\Lambda'^{(t)})\in\calb_{Z^{(t)},Z'^{(t)}}}\rho_{\Lambda^{(t)}}\otimes\rho_{\Lambda'^{(t)}} \\
&=C_t\sum_{(\Lambda^{(t)},\Lambda'^{(t)})\in\calb_{Z^{(t)},Z'^{(t)}}}R_{\Lambda^{(t)}}\otimes R_{\Lambda'^{(t)}} \\
&=\widetilde f_t\otimes\widetilde f'_t\Biggl(\sum_{(\Lambda,\Lambda')\in\calb_{Z,Z'}}R_\Lambda\otimes R_{\Lambda'}\Biggr)
\end{align*}
for some nonzero constant $C_t$.
Because $\widetilde f_t\otimes\widetilde f'_t$ is an isometry,
the proposition is proved.
\end{proof}

\begin{exam}
Consider $Z=\binom{5,3,1}{3,1}\in\cals_{\Sp_{18}}$ and $Z'=\binom{5,3,1}{4,2,0}\in\cals_{\rmO^\pm_{18}}$.
Then $Z_\rmI=\binom{5}{-}$, $\deg(Z)=0$, $|\cals_Z|=1$; $Z'_\rmI=Z'$, $\deg(Z')=3$, $|\cals_{Z'}|=2^6=64$,
and it is not difficult to check that
\[
\calb_{Z,Z'}=\left\{\,(Z,\Lambda')\mid\Lambda'=\textstyle\binom{5,3,1}{4,2,0},\textstyle\binom{5,3,0}{4,2,1},
\textstyle\binom{4,3,1}{5,2,0},\textstyle\binom{4,3,0}{5,2,1},
\textstyle\binom{5,2,1}{4,3,0},\textstyle\binom{5,2,0}{4,3,1},
\textstyle\binom{4,2,1}{5,3,0},\textstyle\binom{4,2,0}{5,3,1}\,\right\},
\]
i.e.,
\begin{align*}
\omega_{Z,Z'}
&=\rho_Z\otimes\left[\rho_{\binom{5,3,1}{4,2,0}}+\rho_{\binom{5,3,0}{4,2,1}}
+\rho_{\binom{4,3,1}{5,2,0}}+\rho_{\binom{4,3,0}{5,2,1}}
+\rho_{\binom{5,2,1}{4,3,0}}+\rho_{\binom{5,2,0}{4,3,1}}
+\rho_{\binom{4,2,1}{5,3,0}}+\rho_{\binom{4,2,0}{5,3,1}}\right] \\
&=\rho_Z\otimes\sum_{M'=\emptyset,\binom{5}{4},\binom{3}{2},\binom{5,3}{4,2}}
\left[\rho_{\Lambda_{M'}}+\rho_{\Lambda_{M'\cup\Psi'}}\right] \\
&=\sqrt 2\left[\rho^{(1)}_Z\otimes\sum_{M'=\emptyset,\binom{5}{4},\binom{3}{2},\binom{5,3}{4,2}}
\rho^{(1)}_{\Lambda_{M'}}\right].
\end{align*}

Now $\Psi_1=\binom{1}{1}$ is an identical pair and $\Psi'_1=\binom{1}{0}$ is in the core of $\calb_{Z,Z'}$.
So we are in Case (I), and then $Z^{(1)}=\binom{4,2}{2}$ and $Z'^{(1)}=\binom{4,2}{3,1}$.
Then $\deg(Z^{(1)})=0$, $|\cals_{Z^{(1)}}|=1$, $\deg(Z'^{(2)})=2$, $|\cals_{Z'^{(1)}}|=2^4=16$, and
\[
\calb_{Z^{(1)},Z'^{(1)}}=\left\{\,(Z^{(1)},\Lambda'^{(1)})\mid\Lambda'^{(1)}=\textstyle\binom{4,2}{3,1},
\textstyle\binom{4,1}{3,2},\textstyle\binom{3,2}{4,1},\textstyle\binom{3,1}{4,2}\,\right\}.
\]
Then we have
\[
(\widetilde f_1\otimes\widetilde f'_1)(\omega_{Z,Z'})
= \sqrt 2\,\rho_{Z^{(1)}}\otimes\left[\rho_{\binom{4,2}{3,1}}+\rho_{\binom{4,1}{3,2}}
+\rho_{\binom{3,2}{4,1}}+\rho_{\binom{3,1}{4,2}}\right]
= 2\left[\rho^{(2)}_{Z^{(1)}}\otimes\sum_{M''=\emptyset,\binom{4}{3}}\rho^{(2)}_{\Lambda_{M''}} \right].
\]

Then $\Psi_2=\binom{2}{2}$ is an identical pair and $\Psi'_2=\binom{2}{1}$ is in the core of $\calb_{Z^{(1)},Z'^{(1)}}$.
Again we are in Case (I), and then $Z^{(2)}=\binom{3}{-}$ and $Z'^{(2)}=\binom{3}{2}$.
Then $\deg(Z^{(2)})=0$, $|\cals_{Z^{(2)}}|=1$, $\deg(Z'^{(2)})=1$, $|\cals_{Z'^{(2)}}|=2^2=4$, and
\[
\calb_{Z^{(2)},Z'^{(2)}}=\left\{(Z^{(2)},\textstyle\binom{3}{2}),(Z^{(2)},\textstyle\binom{2}{3})\right\}.
\]
Then we have
\[
(\widetilde f_2\otimes\widetilde f'_2)(\omega_{Z,Z'})
= 2\,\rho_{Z^{(2)}}\otimes\left[\rho_{\binom{3}{2}}+\rho_{\binom{2}{3}}\right]
= 2\sqrt 2\,\left[\rho^{(3)}_{Z^{(2)}}\otimes\rho^{(3)}_{Z'^{(2)}} \right].
\]

Then $\Psi_3=\emptyset$, and $\Psi'_3=\binom{3}{2}$ is in the core of $\calb_{Z^{(2)},Z'^{(2)}}$.
Now we are in Case (III), and then $Z^{(3)}=\binom{3}{-}$ and $Z'^{(3)}=\binom{-}{-}$.
Then $\deg(Z^{(3)})=\deg(Z'^{(3)})=0$, $|\cals_{Z^{(3)}}|=|\cals_{Z'^{(3)}}|=1$, and
\[
\calb_{Z^{(3)},Z'^{(3)}}=\{(Z^{(3)},Z'^{(3)})\}.
\]
Now $\calb_{Z^{(3)},Z'^{(3)}}$ is regular, and by Lemma~\ref{0303} we have
\[
(\widetilde f_3\otimes\widetilde f'_3)(\omega_{Z,Z'})
= 2\sqrt 2\,\left[\rho_{Z^{(3)}}\otimes\rho_{Z'^{(3)}} \right]
= 2\sqrt 2\,\left[R_{Z^{(3)}}\otimes R_{Z'^{(3)}} \right].
\]
Therefore, we conclude that
\[
\omega_{Z,Z'}
=R_Z\otimes\left[R_{\binom{5,3,1}{4,2,0}}+R_{\binom{5,3,0}{4,2,1}}
+R_{\binom{4,3,1}{5,2,0}}+R_{\binom{4,3,0}{5,2,1}}
+R_{\binom{5,2,1}{4,3,0}}+R_{\binom{5,2,0}{4,3,1}}
+R_{\binom{4,2,1}{5,3,0}}+R_{\binom{4,2,0}{5,3,1}}\right].
\]
\end{exam}

\section{Lusztig Correspondence on Almost Characters}

\subsection{Lusztig correspondence on almost characters}
If $\bfG=\Sp_{2n},\rmO^\epsilon_{2n},\SO_{2n+1}$ (where $\epsilon=+$ or $-$)
and $s$ is a rational semisimple element in the dual group of $\bfG$,
then groups $\bfG^{(0)}(s),\bfG^{(-)}(s),\bfG^{(+)}(s)$ are defined in \cite{pan-Lusztig-correspondence}
\S 2.2 as follows:

\begin{itemize}
\item If $\bfG=\rmO^\epsilon_{2n}$, then
\begin{itemize}
\item $\bfG^{(0)}(s)$ is a product of general linear groups or unitary groups;

\item $\bfG^{(-)}(s)=\rmO^{\epsilon^{(-)}}_{2n^{(-)}}$;

\item $\bfG^{(+)}(s)=\rmO^{\epsilon^{(+)}}_{2n^{(+)}}$
\end{itemize}
where $\epsilon^{(-)}\epsilon^{(+)}=\epsilon$.

\item If $\bfG=\Sp_{2n}$, then
\begin{itemize}
\item $\bfG^{(0)}(s)$ is a product of general linear groups or unitary groups;

\item $\bfG^{(-)}(s)=\rmO^{\epsilon^{(-)}}_{2n^{(-)}}$;

\item $\bfG^{(+)}(s)=\Sp_{2n^{(+)}}$.
\end{itemize}

\item If $\bfG=\SO_{2n+1}$, then
\begin{itemize}
\item $\bfG^{(0)}(s)$ is a product of general linear groups or unitary groups;

\item $\bfG^{(-)}(s)=\Sp_{2n^{(-)}}$;

\item $\bfG^{(+)}(s)=\Sp_{2n^{(+)}}$.
\end{itemize}
\end{itemize}
Here $n^{(-)},n^{(+)}$ denote the multiplicities of ``eigenvalues'' $-1,+1$ of $s$ respectively.
Then $s$ can be written as
\begin{equation}
s=s^{(0)}\times s^{(-)}\times s^{(+)}
\end{equation}
where $s^{(-)}$ (resp.~$s^{(+)}$) is the part whose eigenvalues are all equal to $-1$ (resp.~$1$),
and $s^{(0)}$ is the part whose eigenvalues do not contain $-1$ or $1$.

We have the \emph{modified Lusztig correspondence} (\cf.~\cite{pan-ambiguity} \S 7, \S 8, \S 9)
\begin{equation}\label{0407}
\Xi_s\colon\cale(\bfG)_s\rightarrow \cale(\bfG^{(0)}(s)\times\bfG^{(-)}(s)\times\bfG^{(+)}(s))_1
\end{equation}
which is bijective and can be extended linearly to be an isometry of inner product spaces
\[
\Xi_s\colon\calv(\bfG)_s\rightarrow \calv(\bfG^{(0)}(s))_1\otimes\calv(\bfG^{(-)}(s))_1\otimes\calv(\bfG^{(+)}(s))_1.
\]
Now we enlarge the isometry $\Xi_s$ to $\Xi_{[s]}$ as follows:
\begin{enumerate}
\item Suppose that $\bfG=\rmO^\epsilon_{2n}$ (where $\epsilon=+$ or $-$)
and a semisimple $s\in\rmO^\epsilon_{2n}(q)$ with
$\bfG^{(-)}(s)=\rmO^{\epsilon^{(-)}}_{2n^{(-)}}$ and $\bfG^{(+)}(s)=\rmO^{\epsilon^{(+)}}_{2n^{(+)}}$ as above.
Let $s'\in\rmO^\epsilon_{2n}(q)$ and $s'',s'''\in\rmO^{-\epsilon}_{2n}(q)$ be the semisimple elements such that
\begin{itemize}
\item $\bfG^{(0)}(s')=\bfG^{(0)}(s)$, $\bfG^{(-)}(s')=\rmO^{-\epsilon^{(-)}}_{2n^{(-)}}$
and $\bfG^{(+)}(s')=\rmO^{-\epsilon^{(+)}}_{2n^{(+)}}$;

\item $\bfG^{(0)}(s'')=\bfG^{(0)}(s)$, $\bfG^{(-)}(s'')=\rmO^{\epsilon^{(-)}}_{2n^{(-)}}$
and $\bfG^{(+)}(s'')=\rmO^{-\epsilon^{(+)}}_{2n^{(+)}}$;

\item $\bfG^{(0)}(s''')=\bfG^{(0)}(s)$, $\bfG^{(-)}(s''')=\rmO^{-\epsilon^{(-)}}_{2n^{(-)}}$
and $\bfG^{(+)}(s''')=\rmO^{\epsilon^{(+)}}_{2n^{(+)}}$.
\end{itemize}
Then we have a bijection
\begin{multline}\label{0401}
\Xi_{[s]}\colon\cale(\rmO^\epsilon_{2n})_s\cup\cale(\rmO^\epsilon_{2n})_{s'}
\cup\cale(\rmO^{-\epsilon}_{2n})_{s''}\cup\cale(\rmO^{-\epsilon}_{2n})_{s'''}
\rightarrow \cale(\bfG^{(0)}(s)\times\rmO^\pm_{2n^{(-)}}\times\rmO^\pm_{2n^{(+)}})_1.
\end{multline}

\item Suppose that $\bfG=\Sp_{2n}$ and a semisimple $s\in\SO_{2n+1}(q)$ with
$\bfG^{(-)}(s)=\rmO^{\epsilon^{(-)}}_{2n^{(-)}}$ as above.
Let $s'\in\SO_{2n+1}(q)$ be a semisimple element such that $\bfG^{(0)}(s')=\bfG^{(0)}(s)$,
$\bfG^{(-)}(s')=\rmO^{-\epsilon^{(-)}}_{2n^{(-)}}$ and $\bfG^{(+)}(s')=\bfG^{(+)}(s)$.
Then we have a bijection
\begin{equation}\label{0402}
\Xi_{[s]}\colon\cale(\Sp_{2n})_s\cup\cale(\Sp_{2n})_{s'}
\rightarrow\cale(\bfG^{(0)}(s)\times\rmO^\pm_{2n^{(-)}}\times\Sp_{2n^{(+)}})_1.
\end{equation}

\item Suppose that $\bfG=\SO_{2n+1}$ and a semisimple $s\in\Sp_{2n}(q)$.
Then we have
\begin{equation}\label{0403}
\Xi_{[s]}=\Xi_s\colon\cale(\SO_{2n+1})_s\rightarrow \calv(\bfG^{(0)}(s)\times\Sp_{2n^{(-)}}\times\Sp_{2n^{(+)}})_1.
\end{equation}
\end{enumerate}

We denote
\begin{align*}
\cale(\bfG)_{[s]} &=\begin{cases}
\cale(\rmO^\epsilon_{2n})_s\cup\cale(\rmO^\epsilon_{2n})_{s'}
\cup\cale(\rmO^{-\epsilon}_{2n})_{s''}\cup\cale(\rmO^{-\epsilon}_{2n})_{s'''}, & \text{if $\bfG=\rmO^\pm_{2n}$};\\
\cale(\bfG)_s\cup\cale(\bfG)_{s'}, & \text{if $\bfG=\Sp_{2n}$};\\
\cale(\bfG)_s, & \text{if $\bfG=\SO_{2n+1}$},
\end{cases} \\
X(\bfG)_{[s]} &=\{\,x\in X(\bfG)\mid\rho_x\in\cale(\bfG)_{[s]}\,\}, \\
\calv(\bfG)_{[s]} &=\text{the subspace spanned by }\cale(\bfG)_{[s]}.
\end{align*}
Now for $\bfG=\rmO^\pm_{2n}$, $\Sp_{2n}$, $\SO_{2n+1}$ we have the groups
$\bfG^{(0)}[s],\bfG^{(-)}[s],\bfG^{(+)}[s]$ as follows.
\begin{itemize}
\item If $\bfG=\rmO^\pm_{2n}$, then
\begin{itemize}
\item $\bfG^{(0)}[s]$ is a product of general linear groups or unitary groups;

\item $\bfG^{(-)}[s]=\rmO^\pm_{2n^{(-)}}$;

\item $\bfG^{(+)}[s]=\rmO^\pm_{2n^{(+)}}$.
\end{itemize}

\item If $\bfG=\Sp_{2n}$, then
\begin{itemize}
\item $\bfG^{(0)}[s]$ is a product of general linear groups or unitary groups;

\item $\bfG^{(-)}[s]=\rmO^\pm_{2n^{(-)}}$;

\item $\bfG^{(+)}[s]=\Sp_{2n^{(+)}}$.
\end{itemize}

\item If $\bfG=\SO_{2n+1}$, then
\begin{itemize}
\item $\bfG^{(0)}[s]$ is a product of general linear groups or unitary groups;

\item $\bfG^{(-)}[s]=\Sp_{2n^{(-)}}$;

\item $\bfG^{(+)}[s]=\Sp_{2n^{(+)}}$.
\end{itemize}
\end{itemize}
Then the bijection (\ref{0407}) is enlarged to be a bijection
\begin{equation}\label{0408}
\Xi_{[s]}\colon\cale(\bfG)_{[s]}\rightarrow \cale(\bfG^{(0)}[s]\times\bfG^{(-)}[s]\times\bfG^{(+)}[s])_1
\end{equation}
still called a (modified) \emph{Lusztig correspondence}, which can be extended to an isometry
\begin{equation}\label{0406}
\Xi_{[s]}\colon\calv(\bfG)_{[s]}\rightarrow\calv(\bfG^{(0)}[s])_1\otimes\calv(\bfG^{(-)}[s])_1\otimes\calv(\bfG^{(+)}[s])_1.
\end{equation}

Via isometry (\ref{0406}) the structure of unipotent almost characters of $\bfG^{(0)}[s]\times\bfG^{(-)}[s]\times\bfG^{(+)}[s]$
induces the structure of almost characters in $\calv(\bfG)_{[s]}$ as follows.
For $x\in X(\bfG)_{[s]}$ with
$\Xi_{[s]}(\rho_x)=\rho_{x^{(0)}}\otimes\rho_{\Lambda_x^{(-)}}\otimes\rho_{\Lambda_x^{(+)}}$,
the almost character $R_x\in\calv(\bfG)_{[s]}$ is defined by
\begin{equation}\label{0421}
R_x=\sum_{\Lambda_{x_1}^{(-)}\in\cals_{\bfG^{(-)}[s]}}\,\sum_{\Lambda_{x_1}^{(+)}\in\cals_{\bfG^{(+)}[s]}}
\langle R_{\Lambda_x^{(-)}},\rho_{\Lambda_{x_1}^{(-)}}\rangle_{\bfG^{(-)}[s]}
\langle R_{\Lambda_x^{(+)}},\rho_{\Lambda_{x_1}^{(+)}}\rangle_{\bfG^{(+)}[s]}\rho_{x_1},
\end{equation}
i.e.,
\[
\langle R_x,\rho_{x_1}\rangle_\bfG=\begin{cases}
\langle R_{\Lambda_x^{(-)}},\rho_{\Lambda_{x_1}^{(-)}}\rangle_{\bfG^{(-)}[s]}
\langle R_{\Lambda_x^{(+)}},\rho_{\Lambda_{x_1}^{(+)}}\rangle_{\bfG^{(+)}[s]}, & \text{if $x,x_1\in X(\bfG)_{[s]}$ for some $s$};\\
0, & \text{otherwise}.
\end{cases}
\]
From this definition, we conclude that
\begin{equation}\label{0414}
\Xi_{[s]}(R_x)=R_{x^{(0)}}\otimes R_{\Lambda_x^{(-)}}\otimes R_{\Lambda_x^{(+)}}
\end{equation}
for $x\in X(\bfG)_{[s]}$, and the mapping of almost characters
\begin{equation}
\Xi_{[s]}\colon\cala(\bfG)_{[s]}\rightarrow \cala(\bfG^{(0)}[s]\times\bfG^{(-)}[s]\times\bfG^{(+)}[s])_1
\end{equation}
is a bijection.
The following corollary follows from (\ref{0421}) immediately.
\begin{cor}
Let $\bfG=\Sp_{2n},\rmO^\pm_{2n},\SO_{2n+1}$.
Then $\rho_x=R_x$ if and only if $\rho_{\Lambda_x^{(-)}}=R_{\Lambda_x^{(-)}}$ and
$\rho_{\Lambda_x^{(+)}}=R_{\Lambda_x^{(+)}}$.
\end{cor}

\begin{exam}\label{0411}
The symplectic group $\Sp_{2n}(q)$, $n\geq1$,
has two irreducible characters $\xi_1,\xi_2$ of degree $\frac{1}{2}(q^n+1)$ and
two irreducible characters $\xi_3,\xi_4$ of degree $\frac{1}{2}(q^n-1)$.
It is known that $\xi_1,\xi_2$ are in $\cale(\Sp_{2n})_s$ for a semisimple element $s\in\SO_{2n+1}(q)$ such that
$\bfG^{(0)}(s)=\rmU_0$, $\bfG^{(-)}(s)=\rmO^+_{2n}$ and $\bfG^{(+)}(s)=\Sp_0$ and
\[
\Xi_s\colon\{\xi_1,\xi_2\}
\rightarrow\bigl\{{\bf1}\otimes{\bf 1}_{\rmO^+_{2n}}\otimes{\bf1},{\bf1}\otimes\sgn_{\rmO^+_{2n}}\otimes{\bf1}\bigr\}
=\bigl\{{\bf1}\otimes\rho_{\binom{n}{0}}\otimes{\bf1},{\bf1}\otimes\rho_{\binom{0}{n}}\otimes{\bf1}\bigr\};
\]
and $\xi_3,\xi_4$ are in $\cale(\Sp_{2n})_{s'}$ for a semisimple element $s'\in\SO_{2n+1}(q)$ such that
$\bfG^{(0)}(s')=\rmU_0$, $\bfG^{(-)}(s')=\rmO^-_{2n}$ and $\bfG^{(+)}(s')=\Sp_0$ and
\[
\Xi_{s'}\colon\{\xi_3,\xi_4\}
\rightarrow\bigl\{{\bf1}\otimes{\bf 1}_{\rmO^-_{2n}}\otimes{\bf1},{\bf1}\otimes\sgn_{\rmO^-_{2n}}\otimes{\bf1}\bigr\}
=\bigl\{{\bf1}\otimes\rho_{\binom{-}{n,0}}\otimes{\bf1},{\bf1}\otimes\rho_{\binom{n,0}{-}}\otimes{\bf1}\bigr\}.
\]
Now the four irreducible unipotent characters
\[
\bigl\{{\bf1}\otimes{\bf 1}_{\rmO^+_{2n}}\otimes{\bf1},{\bf1}\otimes\sgn_{\rmO^+_{2n}}\otimes{\bf1},
{\bf1}\otimes{\bf 1}_{\rmO^-_{2n}}\otimes{\bf1},{\bf1}\otimes\sgn_{\rmO^-_{2n}}\otimes{\bf1}\bigr\}
\]
of $\rmU_0(q)\times\rmO^\pm_{2n}(q)\times\Sp_0(q)$ form a family.
So by the modified Lusztig correspondence,
the four irreducible characters $\{\xi_1,\xi_2,\xi_3,\xi_4\}$ of $\Sp_{2n}(q)$ also form a family.
If we fix choices of $\Xi_s$ such that $\Xi_s(\xi_1)={\bf1}\otimes{\bf 1}_{\rmO^+_{2n}}\otimes{\bf1}$ and
$\Xi_{s'}(\xi_3)={\bf1}\otimes{\bf 1}_{\rmO^-_{2n}}\otimes{\bf1}$.
Then we have four almost characters $\{R_1,R_2,R_3,R_4\}$ of the family such that
\begin{align*}
R_1 &= \frac{1}{2}(\xi_1+\xi_2+\xi_3+\xi_4), & R_2 &= \frac{1}{2}(\xi_1+\xi_2-\xi_3-\xi_4), \\
R_3 &= \frac{1}{2}(\xi_1-\xi_2+\xi_3-\xi_4), & R_4 &= \frac{1}{2}(\xi_1-\xi_2-\xi_3+\xi_4)
\end{align*}
and then
\begin{align*}
\Xi_{[s]}(R_1) &= {\bf 1}\otimes R_{\binom{n}{0}}\otimes{\bf1},
& \Xi_{[s]}(R_2) &= {\bf 1}\otimes R_{\binom{0}{n}}\otimes{\bf1}, \\
\Xi_{[s]}(R_3) &= {\bf 1}\otimes R_{\binom{-}{n,0}}\otimes{\bf1},
& \Xi_{[s]}(R_4) &= {\bf 1}\otimes R_{\binom{n,0}{-}}\otimes{\bf1}.
\end{align*}
Note that $R_1+R_3=\xi_1+\xi_3$.
\end{exam}

\subsection{Regular character $\grR_\bfG$}
For $\bfG=\Sp_{2n}$ or $\rmO^\pm_{2n}$, we define
\[
\grR_\bfG=\sum_{x\in X(\bfG)}\rho_x\otimes\rho_x\in\calv(\bfG\times\bfG).
\]
Note that $\grR_{\rmO^\pm_{2n}}\in\calv(\rmO_{2n}^+\times\rmO^+_{2n})\oplus\calv(\rmO_{2n}^-\times\rmO^-_{2n})$.
Let $\grR_{\bfG,1}$ denote the unipotent part of $\grR_\bfG$, i.e.,
\[
\grR_{\bfG,1}=\sum_{\Lambda\in\cals_\bfG}\rho_\Lambda\otimes\rho_\Lambda\in\calv(\bfG\times\bfG)_1.
\]

\begin{lem}\label{0404}
Let $Z$ be a special symbol of rank $n$ and of defect $1$ or $0$.
Then
\[
\sum_{\Lambda\in\cals_Z}\rho_\Lambda\otimes\rho_\Lambda
=\sum_{\Lambda\in\cals_Z}R_\Lambda\otimes R_\Lambda.
\]
\end{lem}
\begin{proof}
Write
\[
\sum_{\Lambda\in\cals_Z}\rho_\Lambda\otimes\rho_\Lambda
=\sum_{\Lambda_1,\Lambda_2\in\cals_Z}m_{\Lambda_1,\Lambda_2} R_{\Lambda_1}\otimes R_{\Lambda_2}
\]
for some $m_{\Lambda_1,\Lambda_2}\in\bbC$.
Then from the proof of Lemma~\ref{0207} and the proof of Lemma~\ref{0208}, we have

\begin{align*}
m_{\Lambda_1,\Lambda_2}
=\biggl\langle\sum_{\Lambda\in\cals_Z}\rho_\Lambda\otimes\rho_\Lambda,
R_{\Lambda_1}\otimes R_{\Lambda_2}\biggr\rangle_{\!\bfG}
&= \frac{1}{2^{2\deg(Z)}}\sum_{\Lambda\in\cals_Z}
(-1)^{\langle\Lambda,\Lambda_1\rangle+\langle\Lambda,\Lambda_2\rangle} \\
&= \begin{cases}
1, & \text{if $\Lambda_1=\Lambda_2$};\\
0, & \text{if $\Lambda_1\neq\Lambda_2$}.
\end{cases}
\end{align*}
Hence the lemma is proved.
\end{proof}

\begin{lem}\label{0405}
For $\bfG=\Sp_{2n}$ or $\rmO^\pm_{2n}$, we have
$\grR_{\bfG,1}=\sum_{\Lambda\in\cals_\bfG}R_\Lambda\otimes R_\Lambda$.
\end{lem}
\begin{proof}
By Lemma~\ref{0404}, we have
\[
\grR_{\bfG,1}
=\sum_{\Lambda\in\cals_\bfG}\rho_\Lambda\otimes\rho_\Lambda
=\sum_Z\sum_{\Lambda\in\cals_Z}\rho_\Lambda\otimes\rho_\Lambda
=\sum_Z\sum_{\Lambda\in\cals_Z}R_\Lambda\otimes R_\Lambda
=\sum_{\Lambda\in\cals_\bfG}R_\Lambda\otimes R_\Lambda.
\]
Here $\sum_Z$ sums over special symbols of rank $n$ and defect $1$ (resp.~defect $0$) if $\bfG=\Sp_{2n}$
(resp.~$\bfG=\rmO^\pm_{2n}$).
\end{proof}

\begin{prop}
For $\bfG=\Sp_{2n}$ or $\rmO^\pm_{2n}$, we have
$\grR_\bfG=\sum_{x\in X(\bfG)}R_x\otimes R_x$.
\end{prop}
\begin{proof}
We can write $\grR_\bfG=\sum_{[s]}\sum_{x\in X(\bfG)_{[s]}}\rho_x\otimes\rho_x$.
Now the mapping
\[
\Xi_{[s]}\colon\calv(\bfG)_{[s]}\rightarrow \calv(\bfG^{(0)}[s]\times\bfG^{(-)}[s]\times\bfG^{(+)}[s])_1
\]
gives
\begin{align*}
& \Xi_{[s]}\Biggl(\sum_{x\in X(\bfG)_{[s]}}\rho_x\otimes\rho_x\Biggr) \\
&=\Biggl[\sum_{x^{(0)}\in X(\bfG^{(0)})_1}\rho_{x^{(0)}}\otimes\rho_{x^{(0)}}\Biggr]
\otimes\Biggl[\sum_{\Lambda_x^{(-)}\in\cals_{\bfG^{(-)}[s]}}\rho_{\Lambda_x^{(-)}}\otimes\rho_{\Lambda_x^{(-)}}\Biggr]
\otimes\Biggl[\sum_{\Lambda_x^{(+)}\in\cals_{\bfG^{(+)}[s]}}\rho_{\Lambda_x^{(+)}}\otimes\rho_{\Lambda_x^{(+)}}\Biggr] \\
&=\Biggl[\sum_{x^{(0)}\in X(\bfG^{(0)})_1}R_{x^{(0)}}\otimes R_{x^{(0)}}\Biggr]
\otimes\Biggl[\sum_{\Lambda_x^{(-)}\in\cals_{\bfG^{(-)}[s]}}R_{\Lambda_x^{(-)}}\otimes R_{\Lambda_x^{(-)}}\Biggr]
\otimes\Biggl[\sum_{\Lambda_x^{(+)}\in\cals_{\bfG^{(+)}[s]}}R_{\Lambda_x^{(+)}}\otimes R_{\Lambda_x^{(+)}}\Biggr] \\
&= \Xi_{[s]}\Biggl(\sum_{x\in X(\bfG)_{[s]}}R_x\otimes R_x\Biggr).
\end{align*}
Because $\Xi_{[s]}$ is an isometry, we have
\[
\sum_{x\in X(\bfG)_{[s]}}\rho_x\otimes\rho_x=\sum_{x\in X(\bfG)_{[s]}}R_x\otimes R_x.
\]
Finally
\[
\sum_{x\in X(\bfG)}\rho_x\otimes\rho_x
=\sum_{[s]}\sum_{x\in X(\bfG)_{[s]}}\rho_x\otimes\rho_x
=\sum_{[s]}\sum_{x\in X(\bfG)_{[s]}}R_x\otimes R_x
=\sum_{x\in X(\bfG)}R_x\otimes R_x.
\]
\end{proof}

\subsection{Lusztig correspondence and finite theta correspondence I}\label{0425}
From \cite{pan-Lusztig-correspondence} \S 6, for $\epsilon=+$ or $-$, we have
\begin{equation}\label{0412}
\omega^\psi_{\Sp_{2n},\rmO^\epsilon_{2n'}}
=\sum_{l=0}^{\min(n,n')}\sum_{(t)\subset X_{l,l}}\omega^\psi_{\Sp_{2n},\rmO^\epsilon_{2n'},t}
\end{equation}
where $\omega^\psi_{\Sp_{2n},\rmO^\epsilon_{2n'},t}$ is the orthogonal projection of
$\omega^\psi_{\Sp_{2n},\rmO^\epsilon_{2n'}}$ on the subspace
$\calv(\Sp_{2n})_{t_n}\otimes\calv(\rmO^\epsilon_{2n'})_{t_{n'}}$,
$t_n=(t,1,\ldots,1)\in\SO_{2n+1}(q)$ ($(n-l)$ copies of $1$),
and $t_{n'}=(t,1,\ldots,1)\in\rmO^\epsilon_{2n'}(q)$ ($(n'-l)$ copies of $1$).
We know that
\begin{itemize}
\item $\bfG^{(0)}(t_n)=\bfG'^{(0)}(t_{n'})=\bfG^{(0)}$ is a product of general linear groups or unitary groups;

\item $\bfG^{(-)}(t_n)=\bfG'^{(-)}(t_{n'})=\rmO^{\epsilon^{(-)}}_{2n^{(-)}}$;

\item $(\bfG^{(+)}(t_n),\bfG'^{(+)}(t_{n'}))=(\Sp_{2n^{(+)}},\rmO^{\epsilon^{(+)}}_{2n'^{(+)}})$
\end{itemize}
for some $\epsilon^{(-)},\epsilon^{(+)}$ such that $\epsilon^{(-)}\epsilon^{(+)}=\epsilon$.
Moreover, from \cite{pan-Lusztig-correspondence} proposition~6.8, we have
\[
\Xi_{t_n}\otimes\Xi_{t_{n'}}\colon\omega^{\psi\sharp}_{\Sp_{2n},\rmO^\epsilon_{2n'},t}
\mapsto \grR^\sharp_{\bfG^{(0)},1}\otimes\grR^\sharp_{\rmO^{\epsilon^{(-)}}_{2n^{(-)}},1}
\otimes\omega^\sharp_{\Sp_{2n^{(+)}},\rmO^{\epsilon^{(+)}}_{2n'^{(+)}},1}
\]
where $f^\sharp$ denotes the uniform projection of a class function $f$ (\cf.~\cite{pan-uniform}).
From \cite{pan-ambiguity} theorem~9.12,
we know that the modified Lusztig correspondence $\Xi_{t_n}\otimes\Xi_{t_{n'}}$ 
can be properly chosen so that the uniform projection in the mapping can be removed.
Then we have 
\[
\Xi_{t_n}\otimes\Xi_{t_{n'}}\colon\omega^\psi_{\Sp_{2n},\rmO^\epsilon_{2n'},t}
\mapsto \grR_{\bfG^{(0)},1}\otimes\grR_{\rmO^{\epsilon^{(-)}}_{2n^{(-)}},1}
\otimes\omega_{\Sp_{2n^{(+)}},\rmO^{\epsilon^{(+)}}_{2n'^{(+)}},1}.
\]

Now there are element $t'\in X_{l,l}$ for the pair $(\Sp_{2n},\rmO^\epsilon_{2n'})$,
and elements $t'',t'''\in X_{l,l}$ for the pair $(\Sp_{2n},\rmO^{-\epsilon}_{2n'})$ such that
\begin{itemize}
\item $\bfG^{(0)}(t'_n)=\bfG^{(0)}(t''_n)=\bfG^{(0)}(t'''_n)=\bfG^{(0)}$;

\item $\bfG^{(-)}(t'_n)=\bfG^{(-)}(t'''_n)=\rmO^{-\epsilon^{(-)}}_{2n^{(-)}}$,
$\bfG^{(-)}(t''_n)=\rmO^{\epsilon^{(-)}}_{2n^{(-)}}$;

\item $\bfG'^{(+)}(t'_{n'})=\bfG'^{(+)}(t''_{n'})=\rmO^{-\epsilon^{(+)}}_{2n'^{(+)}}$,
$\bfG'^{(+)}(t'''_{n'})=\rmO^{\epsilon^{(+)}}_{2n'^{(+)}}$.
\end{itemize}
Then we have 
\begin{align*}
\Xi_{t'_n}\otimes\Xi_{t'_{n'}}\colon\omega^\psi_{\Sp_{2n},\rmO^\epsilon_{2n'},t'}
&\mapsto \grR_{\bfG^{(0)},1}\otimes\grR_{\rmO^{-\epsilon^{(-)}}_{2n^{(-)}},1}
\otimes\omega_{\Sp_{2n^{(+)}},\rmO^{-\epsilon^{(+)}}_{2n'^{(+)}},1}; \\
\Xi_{t''_n}\otimes\Xi_{t''_{n'}}\colon\omega^\psi_{\Sp_{2n},\rmO^{-\epsilon}_{2n'},t''}
&\mapsto \grR_{\bfG^{(0)},1}\otimes\grR_{\rmO^{\epsilon^{(-)}}_{2n^{(-)}},1}
\otimes\omega_{\Sp_{2n^{(+)}},\rmO^{-\epsilon^{(+)}}_{2n'^{(+)}},1}; \\
\Xi_{t'''_n}\otimes\Xi_{t'''_{n'}}\colon\omega^\psi_{\Sp_{2n},\rmO^{-\epsilon}_{2n'},t'''}
&\mapsto \grR_{\bfG^{(0)},1}\otimes\grR_{\rmO^{-\epsilon^{(-)}}_{2n^{(-)}},1}
\otimes\omega_{\Sp_{2n^{(+)}},\rmO^{\epsilon^{(+)}}_{2n'^{(+)}},1}.
\end{align*}
Let $[t]$ denote the union of the Weyl group orbits $(t),(t'),(t''),(t''')$ in $X_{l,l}$ and
\[
\omega^\psi_{\Sp_{2n},\rmO^\pm_{2n'},[t]}
=\omega^\psi_{\Sp_{2n},\rmO^\epsilon_{2n'},t}
+\omega^\psi_{\Sp_{2n},\rmO^\epsilon_{2n'},t'}
+\omega^\psi_{\Sp_{2n},\rmO^{-\epsilon}_{2n'},t''}
+\omega^\psi_{\Sp_{2n},\rmO^{-\epsilon}_{2n'},t'''}.
\]
Then from (\ref{0412}), we have
\begin{equation}\label{0415}
\omega^\psi_{\Sp_{2n},\rmO^\pm_{2n'}}
=\sum_{l=0}^{\min(n,n')}\sum_{[t]\subset X_{l,l}}\omega^\psi_{\Sp_{2n},\rmO^\pm_{2n'},[t]}
\end{equation}
and the isometry
\begin{multline}\label{0413}
\Xi_{[t_n]}\otimes\Xi_{[t_{n'}]}\colon\calv(\Sp_{2n})_{[t_n]}\otimes\calv(\rmO^\pm_{2n'})_{[t_{n'}]} \\
\rightarrow \calv(\bfG^{(0)}\times\bfG^{(0)})_1\otimes\calv(\rmO^\pm_{2n^{(-)}}\times\rmO^\pm_{2n^{(-)}})_1
\otimes\calv(\Sp_{2n^{(+)}}\times\rmO^\pm_{2n'^{(+)}})_1
\end{multline}
satisfies
\begin{equation}\label{0420}
\omega^\psi_{\Sp_{2n},\rmO^\pm_{2n'},[t]}
\mapsto \grR_{\bfG^{(0)},1}\otimes\grR_{\rmO^\pm_{2n^{(-)}},1}
\otimes\omega_{\Sp_{2n^{(+)}},\rmO^\pm_{2n'^{(+)}},1}.
\end{equation}
\begin{prop}\label{0416}
Let $(\bfG,\bfG')=(\Sp_{2n},\rmO^\pm_{2n'})$ and $t\in X_{l,l}$ for some $l\leq\min(n,n')$.
Then
\[
\omega^\psi_{\bfG,\bfG',[t]}=\sum_{(x,x')}R_x\otimes R_{x'}.
\]
where $\sum_{(x,x')}$ sums over all $(x,x')\in X(\bfG)_{[t_n]}\times X(\bfG')_{[t_{n'}]}$ such that
\begin{itemize}
\item $x^{(0)}=x'^{(0)}\in X(\bfG^{(0)})_1=X(\bfG'^{(0)})_1$;

\item $\Lambda_x^{(-)}=\Lambda_{x'}^{(-)}\in\cals_{\rmO^\pm_{2n^{(-)}}}$;

\item $(\Lambda_x^{(+)},\Lambda_{x'}^{(+)})\in\calb_{\Sp_{2n^{(+)}},\rmO^\pm_{2n'^{(+)}}}$.
\end{itemize}
\end{prop}
\begin{proof}
By Lemma~\ref{0405} and Corollary~\ref{0306}, we have
\begin{align*}
& \grR_{\bfG^{(0)},1}\otimes\grR_{\rmO^\pm_{2n^{(-)}},1}\otimes\omega_{\Sp_{2n^{(+)}},\rmO^\pm_{2n'^{(+)}},1} \\
&= \sum_{x^{(0)}\in X(\bfG^{(0)})_1}\,\sum_{\Lambda_x^{(-)}\in\cals_{\rmO^\pm_{2n^{(-)}}}}\,
\sum_{(\Lambda_x^{(+)},\Lambda_{x'}^{(+)})\in\calb_{\Sp_{2n^{(+)}},\rmO^\pm_{2n'^{(+)}}}}\!\!\!\!\!\!\!\!\!\!\!
\bigl[\rho_{x^{(0)}}\otimes\rho_{\Lambda_x^{(-)}}\otimes\rho_{\Lambda_x^{(+)}}\bigr]
\otimes\bigl[\rho_{x^{(0)}}\otimes\rho_{\Lambda_x^{(-)}}\otimes\rho_{\Lambda_{x'}^{(+)}}\bigr] \\
&= \sum_{x^{(0)}\in X(\bfG^{(0)})_1}\,\sum_{\Lambda_x^{(-)}\in\cals_{\rmO^\pm_{2n^{(-)}}}}\,
\sum_{(\Lambda_x^{(+)},\Lambda_{x'}^{(+)})\in\calb_{\Sp_{2n^{(+)}},\rmO^\pm_{2n'^{(+)}}}}\!\!\!\!\!\!\!\!\!\!\!
\bigl[R_{x^{(0)}}\otimes R_{\Lambda_x^{(-)}}\otimes R_{\Lambda_x^{(+)}}\bigr]
\otimes\bigl[R_{x^{(0)}}\otimes R_{\Lambda_x^{(-)}}\otimes R_{\Lambda_{x'}^{(+)}}\bigr].
\end{align*}
Because of (\ref{0413}) and (\ref{0414}),
the proposition is obtained by taking the inverse image $\Xi_{[t_n]}^{-1}\otimes\Xi_{[t_{n'}]}^{-1}$
and by (\ref{0420}).
\end{proof}

From the proposition, we see that $\omega^\psi_{\Sp_{2n},\rmO^\pm_{2n'}}$
is a sum of terms $R_x\otimes R_{x'}$ for $(x,x')$ satisfying certain condition.
We say that $R_x\otimes R_{x'}$ \emph{occurs in} $\omega_{\Sp_{2n},\rmO^\pm_{2n'}}$
if it appears in the summation.

\begin{cor}\label{0409}
Let $(\bfG,\bfG')=(\Sp_{2n},\rmO^\pm_{2n'})$,
$R_x\in\cala(\bfG)_{[s]}$ and $R_{x'}\in\cala(\bfG')_{[s']}$ for some semisimple elements $s,s'$.
Write
\[
\Xi_{[s]}(R_x)=R_{x^{(0)}}\otimes R_{\Lambda_x^{(-)}}\otimes R_{\Lambda_x^{(+)}}\quad\text{and}\quad
\Xi_{[s']}(R_{x'})=R_{x'^{(0)}}\otimes R_{\Lambda_{x'}^{(-)}}\otimes R_{\Lambda_{x'}^{(+)}}.
\]
Then $R_x\otimes R_{x'}$ occurs in $\omega^\psi_{\bfG,\bfG'}$ if and only if
the following conditions hold:
\begin{itemize}
\item $s^{(0)}=s'^{(0)}$ (up to conjugation), 
$\bfG^{(0)}[s]\simeq \bfG'^{(0)}[s']$ and $R_{x^{(0)}}=R_{x'^{(0)}}$;

\item $\bfG^{(-)}[s]\simeq \bfG'^{(-)}[s']$ and $R_{\Lambda_x^{(-)}}=R_{\Lambda_{x'}^{(-)}}$;

\item $R_{\Lambda_x^{(+)}}\otimes R_{\Lambda_{x'}^{(+)}}$ occurs in $\omega_{\bfG^{(+)}[s],\bfG'^{(+)}[s'],1}$.
\end{itemize}
\end{cor}
\begin{proof}
The corollary follows from (\ref{0415}) and Proposition~\ref{0416} immediately.
\end{proof}

\begin{proof}[Proof of Theorem~\ref{0101}: Part I]
Let $(\bfG,\bfG')=(\Sp_{2n},\rmO^\pm_{2n'})$, $\rho_x\in\cale(\bfG)_s$, $\rho_{x'}\in\cale(\bfG')_{s'}$
for some $s,s'$, and write
\[
\Xi_{[s]}(\rho_x)=\rho_{x^{(0)}}\otimes\rho_{\Lambda_x^{(-)}}\otimes\rho_{\Lambda_x^{(+)}}\quad\text{and}\quad
\Xi_{[s']}(\rho_{x'})=\rho_{x'^{(0)}}\otimes\rho_{\Lambda_{x'}^{(-)}}\otimes\rho_{\Lambda_{x'}^{(+)}}.
\]
First suppose that $\rho_x\otimes\rho_{x'}$ occurs in $\omega^\psi_{\bfG,\bfG'}$.
Then by \cite{pan-ambiguity} theorem 9.12,
we know that $s^{(0)}=s'^{(0)}$, $\rho_{x^{(0)}}=\rho_{x'^{(0)}}$,
$\rho_{\Lambda_x^{(-)}}=\rho_{\Lambda_{x'}^{(-)}}$,
and $\rho_{\Lambda_x^{(+)}}\otimes\rho_{\Lambda_{x'}^{(+)}}$ occurs in $\omega_{\bfG^{(+)}[s],\bfG'^{(+)}[s'],1}$.
Then we have $R_{x^{(0)}}=R_{x'^{(0)}}$,
$R_{\Lambda_x^{(-)}}=R_{\Lambda_{x'}^{(-)}}$,
and $R_{\Lambda_x^{(+)}}\otimes R_{\Lambda_{x'}^{(+)}}$ occurs in $\omega_{\bfG^{(+)}[s],\bfG'^{(+)}[s'],1}$
by Corollary~\ref{0306}.
Then by Corollary~\ref{0409}, we conclude that $R_x\otimes R_{x'}$ occurs in $\omega^\psi_{\bfG,\bfG'}$.
On the other hand, if $R_x\otimes R_{x'}$ occurs in $\omega^\psi_{\bfG,\bfG'}$,
then by the same argument, we also conclude that $\rho_x\otimes\rho_{x'}$ occurs in $\omega^\psi_{\bfG,\bfG'}$.
\end{proof}

\subsection{Lusztig correspondence and finite theta correspondence II}
In this subsection, let $(\bfG,\bfG')=(\Sp_{2n},\rmO^\epsilon_{2n'+1})$.
Then from \cite{pan-Lusztig-correspondence} section 7, we have
\begin{equation}\label{0419}
\omega^\psi_{\Sp_{2n},\SO_{2n'+1}}\cdot(1\otimes\chi_{\SO_{2n'+1}})
=\sum_{l=0}^{\min(n,n')}\sum_{(t)\subset X_{l,l}^\flat}
\left(\omega^\psi_{\Sp_{2n},\SO_{2n'+1},t,+}+\omega^\psi_{\Sp_{2n},\SO_{2n'+1},t,-}\right)
\end{equation}
where $\omega^\psi_{\Sp_{2n},\SO_{2n'+1},t,\epsilon}$ denotes the orthogonal projection of
$\omega^\psi_{\Sp_{2n},\SO_{2n'+1}}\cdot(1\otimes\chi_{\SO_{2n'+1}})$ over
$\calv(\Sp_{2n})_{t_{n,\epsilon}^\flat}\otimes\calv(\SO_{2n'+1})_{t_{n'}^\flat}$,
$t^\flat_{n,\epsilon}$ is the semisimple element $(t,-1)\in\SO_{2l+1}\times\SO_{2(n-l)}^\epsilon\subseteq\SO_{2n+1}$,
$t^\flat_{n'}$ is the semisimple element $(t,-1)\in\Sp_{2l}\times\Sp_{2(n'-l)}\subseteq\Sp_{2n'}$,
$\chi_{\SO_{2n'+1}}$ is the character of the spinor norm of $\SO_{2n'+1}$.
We know that
\begin{itemize}
\item $\bfG^{(0)}(t^\flat_{n,\epsilon})=\bfG'^{(0)}(t^\flat_{n'})=\bfG^{(0)}$ is a product of general linear groups or unitary groups;

\item $(\bfG^{(-)}(t^\flat_{n,\epsilon}),\bfG'^{(-)}(t^\flat_{n'}))=(\rmO^{\epsilon}_{2n^{(-)}},\Sp_{2n'^{(-)}})$

\item $\bfG^{(+)}(t^\flat_{n,\epsilon})=\bfG'^{(+)}(t^\flat_{n'})=\Sp_{2n^{(+)}}$.
\end{itemize}

As in Subsection~\ref{0425}, by \cite{pan-Lusztig-correspondence} proposition~7.9 and
\cite{pan-ambiguity} theorem 9.9, 
the modified Lusztig correspondences can be properly chosen so that  
\begin{align*}
\Xi_{t^\flat_{n,+}}\otimes\Xi_{t^\flat_{n'}}\colon\omega^\psi_{\Sp_{2n},\SO_{2n'+1},t,+}
&\mapsto \grR_{\bfG^{(0)},1}\otimes\omega^\psi_{\rmO^+_{2n^{(-)}},\Sp_{2n'^{(-)}},1}
\otimes\grR_{\Sp_{2n^{(+)}},1}, \\
\Xi_{t^\flat_{n,-}}\otimes\Xi_{t^\flat_{n'}}\colon\omega^\psi_{\Sp_{2n},\SO_{2n'+1},t,-}
&\mapsto \grR_{\bfG^{(0)},1}\otimes\omega^\psi_{\rmO^-_{2n^{(-)}},\Sp_{2n'^{(-)}},1}
\otimes\grR_{\Sp_{2n^{(+)}},1}.
\end{align*}
Define
\begin{equation}
\omega^\psi_{\Sp_{2n},\SO_{2n'+1},[t]}
=\omega^\psi_{\Sp_{2n},\SO_{2n'+1},t,+}+\omega^\psi_{\Sp_{2n},\SO_{2n'+1},t,-}.
\end{equation}
Then (\ref{0419}) can be rewritten as
\begin{equation}\label{0418}
\omega^\psi_{\Sp_{2n},\SO_{2n'+1}}\cdot(1\otimes\chi_{\SO_{2n'+1}})
=\sum_{l=0}^{\min(n,n')}\sum_{[t]\subset X^\flat_{l,l}}\omega^\psi_{\Sp_{2n},\SO_{2n'+1},[t]}
\end{equation}
and the isometry
\begin{multline}
\Xi_{[t^\flat_n]}\otimes\Xi_{[t^\flat_{n'}]}\colon \calv(\Sp_{2n})_{[t^\flat_n]}\otimes\calv(\SO_{2n'+1})_{[t^\flat_{n'}]} \\
\rightarrow \calv(\bfG^{(0)}\times\bfG^{(0)})_1\otimes\calv(\rmO^\pm_{2n^{(-)}}\times\Sp_{2n'^{(-)}})_1
\otimes\calv(\Sp_{2n^{(+)}}\times\Sp_{2n^{(+)}})_1
\end{multline}
satisfies
\[
\omega^\psi_{\Sp_{2n},\SO_{2n'+1},[t]}
\mapsto \grR_{\bfG^{(0)},1}\otimes\omega^\psi_{\rmO^\pm_{2n^{(-)}},\Sp_{2n'^{(-)}},1}
\otimes\grR_{\Sp_{2n^{(+)}},1}
\]
where $\calv(\Sp_{2n})_{[t^\flat_n]}:=\calv(\Sp_{2n})_{t^\flat_{n,+}}\oplus \calv(\Sp_{2n})_{t^\flat_{n,-}}$
and $\calv(\SO_{2n'+1})_{[t^\flat_{n'}]}:=\calv(\SO_{2n'+1})_{t^\flat_{n'}}$.

\begin{prop}\label{0417}
Let $(\bfG,\bfG')=(\Sp_{2n},\SO_{2n'+1})$ and $t\in X^\flat_{l,l}$ for some $l\leq\min(n,n')$.
Then
\[
\omega^\psi_{\bfG,\bfG',[t]}=\sum_{(x,x')}R_x\otimes R_{x'}.
\]
where $\sum_{(x,x')}$ sums over all $(x,x')\in X(\bfG)_{[t^\flat_n]}\times X(\bfG')_{[t^\flat_{n'}]}$ such that
\begin{itemize}
\item $x^{(0)}=x'^{(0)}\in X(\bfG^{(0)})_1=X(\bfG'^{(0)})_1$;

\item $(\Lambda_x^{(-)},\Lambda_{x'}^{(-)})\in\calb_{\rmO^\pm_{2n^{(-)}},\Sp_{2n'^{(-)}}}$;

\item $\Lambda_x^{(+)}=\Lambda_{x'}^{(+)}\in\cals_{\Sp_{2n^{(+)}}}$.
\end{itemize}
\end{prop}
\begin{proof}
The proof is the same as that of Proposition~\ref{0416}.
\end{proof}

\begin{lem}\label{0422}
Let $R_x\in\cala(\SO_{2n+1})$.
Then $\chi_{\SO_{2n+1}} R_x=R_{\bar x}$ where $\bar x\in X(\SO_{2n+1})$
such that $\bar x^{(0)}=x^{(0)}$,
$\Lambda^{(-)}_{\bar x}=\Lambda^{(+)}_x$ and $\Lambda^{(+)}_{\bar x}=\Lambda_x^{(-)}$.
\end{lem}
\begin{proof}
Suppose that $\bfG=\SO_{2n+1}$ and $R_x\in\cala(\bfG)_{[s]}$ for some $[s]$.
It is known that for $x_1\in X(\bfG)_{[s]}$,
then $\chi_{\SO_{2n+1}}\rho_{x_1}=\rho_{\bar x_1}\in\cale(\bfG)_{[-s]}$
where $\bar x_1\in X(\bfG)_{[-s]}$ satisfying the requirements of the lemma, i.e.,
$\bar x_1^{(0)}=x_1^{(0)}$,
$\Lambda^{(-)}_{\bar x_1}=\Lambda^{(+)}_{x_1}$ and $\Lambda^{(+)}_{\bar x_1}=\Lambda_{x_1}^{(-)}$.
Note that $\bfG^{(-)}[s]=\bfG^{(+)}[-s]$ and $\bfG^{(+)}[s]=\bfG^{(-)}[-s]$.
Then by (\ref{0421}), we have
\begin{align*}
\chi_{\SO_{2n+1}}R_x
&= \sum_{\Lambda_{x_1}^{(-)}\in\cals_{\bfG^{(-)}[s]}}\,\sum_{\Lambda_{x_1}^{(+)}\in\cals_{\bfG^{(+)}[s]}}
\langle R_{\Lambda_x^{(-)}},\rho_{\Lambda_{x_1}^{(-)}}\rangle_{\bfG^{(-)}[s]}
\langle R_{\Lambda_x^{(+)}},\rho_{\Lambda_{x_1}^{(+)}}\rangle_{\bfG^{(+)}[s]}\chi_{\SO_{2n+1}}\rho_{x_1} \\
&= \sum_{\Lambda_{\bar x_1}^{(+)}\in\cals_{\bfG^{(+)}[-s]}}\,\sum_{\Lambda_{\bar x_1}^{(-)}\in\cals_{\bfG^{(-)}[-s]}}
\langle R_{\Lambda_{\bar x}^{(+)}},\rho_{\Lambda_{\bar x_1}^{(+)}}\rangle_{\bfG^{(+)}[-s]}
\langle R_{\Lambda_{\bar x}^{(-)}},\rho_{\Lambda_{\bar x_1}^{(-)}}\rangle_{\bfG^{(-)}[-s]}
\rho_{\bar x_1} \\
&= R_{\bar x}.
\end{align*}
\end{proof}

\begin{prop}\label{0410}
Let $(\bfG,\bfG')=(\Sp_{2n},\SO_{2n'+1})$,
$R_x\in\cala(\bfG)_{[s]}$ and $R_{x'}\in\cala(\bfG')_{[s']}$ for some semisimple elements $s,s'$.
Write
\[
\Xi_{[s]}(R_x)=R_{x^{(0)}}\otimes R_{\Lambda_x^{(-)}}\otimes R_{\Lambda_x^{(+)}}\quad\text{and}\quad
\Xi_{[s']}(R_{x'})=R_{x'^{(0)}}\otimes R_{\Lambda_{x'}^{(-)}}\otimes R_{\Lambda_{x'}^{(+)}}.
\]
Then $R_x\otimes R_{x'}$ occurs in $\omega^\psi_{\bfG,\bfG'}$ if and only if
the following conditions hold:
\begin{itemize}
\item $s^{(0)}=-s'^{(0)}$ (up to conjugation), 
$\bfG^{(0)}[s]\simeq \bfG'^{(0)}[s']$ and $R_{x^{(0)}}=R_{x'^{(0)}}$;

\item $R_{\Lambda_x^{(-)}}\otimes R_{\Lambda_{x'}^{(+)}}$ occurs in $\omega_{\bfG^{(-)}[s],\bfG'^{(+)}[s'],1}$,

\item $\bfG^{(+)}[s]\simeq \bfG'^{(-)}[s']$ and $R_{\Lambda_x^{(+)}}=R_{\Lambda_{x'}^{(-)}}$.
\end{itemize}
\end{prop}
\begin{proof}
First suppose that $R_x\otimes R_{x'}$ occurs in $\omega^\psi_{\Sp_{2n},\SO_{2n'+1}}$.
Then $R_x\otimes(\chi_{\SO_{2n'+1}} R_{x'})$ occurs in 
$\omega^\psi_{\Sp_{2n},\SO_{2n'+1}}\cdot({\bf 1}\otimes\chi_{\SO_{2n'+1}})$,
and by (\ref{0418}), $R_x\otimes(\chi_{\SO_{2n'+1}} R_{x'})$ occurs in 
$\omega^\psi_{\Sp_{2n},\SO_{2n'+1},[t]}$ for some $t$.
By Lemma~\ref{0422}, we know that $\chi_{\SO_{2n'+1}} R_{x'}=R_{\bar x'}$ where $\bar x'\in X(\bfG')_{[-s']}$.
Now $[s]=[t^\flat_n]$ and $[-s']=[t^\flat_{n'}]$ and by Proposition~\ref{0417} we have
\begin{itemize}
\item $x^{(0)}=\bar x'^{(0)}\in X(\bfG^{(0)})_1=X(\bfG'^{(0)})_1$;

\item $(\Lambda_x^{(-)},\Lambda_{\bar x'}^{(-)})\in\calb_{\bfG^{(-)}[s],\bfG'^{(-)}[-s']}$.

\item $\Lambda_x^{(+)}=\Lambda_{\bar x'}^{(+)}\in\cals_{\bfG^{(+)}[s]}=\cals_{\bfG'^{(+)}[-s']}$;
\end{itemize}
Now $\bar x'^{(0)}=x'^{(0)}$,
$\Lambda^{(-)}_{\bar x'}=\Lambda^{(+)}_{x'}$ and $\Lambda^{(+)}_{\bar x'}=\Lambda_{x'}^{(-)}$,
$\bfG'^{(-)}[-s']=\bfG'^{(+)}[s']$ and $\bfG'^{(+)}[-s']=\bfG'^{(-)}[s']$.
Therefore, by Corollary~\ref{0306}, we conclude that
\begin{itemize}
\item $s^{(0)}=-s'^{(0)}$ (up to conjugation), $\bfG^{(0)}[s]\simeq \bfG'^{(0)}[s']$ and $R_{x^{(0)}}=R_{x'^{(0)}}$;

\item $R_{\Lambda_x^{(-)}}\otimes R_{\Lambda_{x'}^{(+)}}$ occurs in $\omega_{\bfG^{(-)}[s],\bfG'^{(+)}[s'],1}$,

\item $\bfG^{(+)}[s]\simeq \bfG'^{(-)}[s']$ and $R_{\Lambda_x^{(+)}}=R_{\Lambda_{x'}^{(-)}}$.
\end{itemize}
It is easy to see that the implication for the other direction is also true from the above argument.
\end{proof}

\begin{exam}
The correspondence for the dual pair $(\Sp_{2n},\SO_1)$, for $n\geq 1$, can be reduced to the correspondence
\[
[{\bf 1}_{\rmO^+_{2n}}+{\bf1}_{\rmO^-_{2n}}]\otimes{\bf 1}_{\Sp_0}
=[\rho_{\binom{n}{0}}+\rho_{\binom{-}{n,0}}]\otimes\rho_{\binom{0}{-}}
=[R_{\binom{n}{0}}+R_{\binom{-}{n,0}}]\otimes R_{\binom{0}{-}}
\]
for $(\rmO^\pm_{2n},\Sp_0)$.
Therefore
\[
\omega_{\Sp_{2n},\SO_1}=(\xi_1+\xi_3)\otimes{\bf 1}_{\SO_1}=(R_1+R_3)\otimes{\bf 1}_{\SO_1}
\]
in the notation of Example~\ref{0411}.
\end{exam}

\begin{proof}[Proof of Theorem~\ref{0101}: Part II]
Let $(\bfG,\bfG')=(\Sp_{2n},\SO^\pm_{2n'+1})$, $\rho_x\in\cale(\bfG)_s$, $\rho_{x'}\in\cale(\bfG')_{s'}$
for some $s,s'$, and write
\[
\Xi_{[s]}(\rho_x)=\rho_{x^{(0)}}\otimes\rho_{\Lambda_x^{(-)}}\otimes\rho_{\Lambda_x^{(+)}}\quad\text{and}\quad
\Xi_{[s']}(\rho_{x'})=\rho_{x'^{(0)}}\otimes\rho_{\Lambda_{x'}^{(-)}}\otimes\rho_{\Lambda_{x'}^{(+)}}.
\]
First suppose that $\rho_x\otimes\rho_{x'}$ occurs in $\omega^\psi_{\bfG,\bfG'}$.
Then by \cite{pan-ambiguity} theorem 9.9,
we know that $s^{(0)}=-s'^{(0)}$, $\rho_{x^{(0)}}=\rho_{x'^{(0)}}$,
$\rho_{\Lambda_x^{(+)}}=\rho_{\Lambda_{x'}^{(-)}}$,
and $\rho_{\Lambda_x^{(-)}}\otimes\rho_{\Lambda_{x'}^{(+)}}$ occurs in $\omega_{\bfG^{(-)}[s],\bfG'^{(+)}[s'],1}$.
Then we have $R_{x^{(0)}}=R_{x'^{(0)}}$,
$R_{\Lambda_x^{(+)}}=R_{\Lambda_{x'}^{(-)}}$,
and $R_{\Lambda_x^{(-)}}\otimes R_{\Lambda_{x'}^{(+)}}$ occurs in
$\omega_{\bfG^{(-)}[s],\bfG'^{(+)}[s'],1}$ by Corollary~\ref{0306}.
Then by Proposition~\ref{0410},
we conclude that $R_x\otimes R_{x'}$ occurs in $\omega^\psi_{\bfG,\bfG'}$.
On the other hand, if $R_x\otimes R_{x'}$ occurs in $\omega^\psi_{\bfG,\bfG'}$,
then by the same argument, we also conclude that $\rho_x\otimes\rho_{x'}$ occurs in $\omega^\psi_{\bfG,\bfG'}$.
\end{proof}

\subsection{Multiplicity one correspondence}
It is well known that the finite theta correspondence is of multiplicity one for the dual pair
$(\Sp_{2n},\rmO^\epsilon_{2n'})$ or $(\Sp_{2n},\rmO_{2n'+1})$ where $\epsilon=+$ or $-$, i.e.,
if we write
\[
\omega^\psi_{\bfG,\bfG'}
=\sum_{x\in X(\bfG),\ x'\in X(\bfG')} m_{x,x'}\rho_x\otimes\rho_{x'},
\]
then $m_{x,x'}=0,1$ for any $(x,x')\in X(\bfG)\times X(\bfG')$.
It should be known to the experts that the correspondence for the pair $(\Sp_{2n},\SO_{2n'+1})$
is still of multiplicity one.
Now we provide a proof for completion.

\begin{prop}\label{0424}
Suppose that we write
\begin{align*}
\omega^\psi_{\Sp_{2n},\SO_{2n'+1}}
&=\sum_{x\in X(\Sp_{2n}),\ x'\in X(\SO_{2n'+1})} m_{x,x'}\rho_x\otimes\rho'_{x'}.
\end{align*}
Then $m_{x,x'}=0,1$ for any $x\in X(\Sp_{2n})$ and any $x'\in X(\SO_{2n'+1})$.
\end{prop}
\begin{proof}
Let $(\bfG,\bfG')=(\Sp_{2n},\SO_{2n'+1})$.
Suppose that
\begin{itemize}
\item $\rho_i\in\cale(\Sp_{2n})_{s_i}$ for some $s_i$ and $i=1,2$,

\item $\rho'\in\cale(\rmO_{2n'+1})_{s'}$ such that $\sgn\cdot\rho'\neq\rho'$,

\item $(\rho_1,\rho'),(\rho_2,\sgn\cdot\rho')\in \Theta_{\Sp_{2n},\rmO_{2n'+1}}$.
\end{itemize}
Because the correspondence for the pair $(\Sp_{2n},\rmO_{2n'+1})$ is of multiplicity one,
we need to show that $\rho_1\neq\rho_2$.
The fact that $\rho_1\neq\rho_2$ follows form \cite{pan-Lusztig-correspondence} immediately.
\end{proof}

However, the multiplicity one property does not hold for the pair $(\Sp_{2n},\SO^\epsilon_{2n'})$.
Here is an easy counterexample:
\begin{exam}\label{0423}
It is easy to see that
\[
\omega_{\Sp_2,\rmO^+_2,1}=[\rho_{\binom{1}{-}}+\rho_{\binom{1,0}{1}}]\otimes\rho_{\binom{1}{0}}
+\rho_{\binom{1}{-}}\otimes\rho_{\binom{0}{1}}.
\]
Because $\rho_{\binom{1}{0}}|_{\SO^+_2}=\rho_{\binom{0}{1}}|_{\SO^+_2}={\bf 1}_{\SO^+_2}$,
we have
\[
\omega_{\Sp_2,\SO^+_2,1}=2({\bf 1}_{\Sp_2}\otimes{\bf 1}_{\SO^+_2})+{\St}_{\Sp_2}\otimes{\bf 1}_{\SO^+_2}
\]
where $\rho_{\binom{1}{-}}={\bf 1}_{\Sp_2}$, and
$\rho_{\binom{1,0}{1}}={\St}_{\Sp_2}$ denotes the Steinberg characters of $\Sp_2(q)$.
Therefore the correspondence for the pair $(\Sp_2,\SO_2^+)$ is not of multiplicity one.
\end{exam}

\section{An Application}

\subsection{$\Theta$-ranks of almost characters}

For $(\bfG,\bfG')=(\Sp_{2n},\rmO^\pm_{2n'})$ or $(\Sp_{2n},\SO_{2n'+1})$,
Theorem~\ref{0101} tells us that
we can define the finite theta correspondence on almost characters as we did for irreducible
characters:
\[
\widetilde\Theta^\psi_{\bfG,\bfG'}=\{\,(R_x,R_{x'})\in\cala(\bfG)\times\cala(\bfG')\mid m_{x,x'}\neq 0\,\}.
\]
Moreover, we know that $(\rho_x,\rho_{x'})\in\Theta^\psi_{\bfG,\bfG'}$ if and only if 
$(R_x,R_{x'})\in\widetilde\Theta^\psi_{\bfG,\bfG'}$.

Now we can define the $\Theta$-rank of an almost character as we did for an irreducible character as follows.
First we define
\[
\sgn_{\rmO^\pm_{2n}}=\sgn_{\rmO^+_{2n}}+\sgn_{\rmO^-_{2n}}\quad\text{and}\quad
\chi_{\rmO^\pm_{2n}}=\chi_{\rmO^+_{2n}}+\chi_{\rmO^-_{2n}},
\]
which are vectors in $\calv(\rmO^\pm_{2n})$.
Then it is easy to see that $\sgn_{\rmO^\pm_{2n}}\cdot\rho_x=\sgn_{\rmO^\epsilon_{2n}}\cdot\rho_x$ and
$\chi_{\rmO^\pm_{2n}}\rho_x=\chi_{\rmO^\epsilon_{2n}}\rho_x$
if $\rho_x\in\cale(\rmO^\epsilon_{2n})$ for $\epsilon=+$ or $-$.
For $x\in X(\rmO^\pm_{2n})$, define $x^\rmt\in X(\rmO^\pm_{2n})$ such that
$(x^\rmt)^{(0)}=x^{(0)}$, $\Lambda_{x^\rmt}^{(-)}=(\Lambda_x^{(-)})^\rmt$
and $\Lambda_{x^\rmt}^{(+)}=(\Lambda_x^{(+)})^\rmt$, also define $\bar x\in X(\rmO^\pm_{2n})$
such that $\bar x^{(0)}=x^{(0)}$, $\Lambda_{\bar x}^{(-)}=\Lambda_x^{(+)}$
and $\Lambda_{\bar x}^{(+)}=\Lambda_x^{(-)}$.
It is known that $\sgn_{\rmO^\pm_{2n}}\rho_x=\rho_{x^\rmt}$ and
$\chi_{\rmO^\pm_{2n}} R_x=R_{\bar x}$.

\begin{lem}
Suppose that $R_x\in\cala(\rmO^\pm_{2n})$.
Then
\[
\sgn_{\rmO^\pm_{2n}}\cdot R_x=R_{x^\rmt}\quad\text{and}\quad
\chi_{\rmO^\pm_{2n}} R_x=R_{\bar x}.
\]
In particular, both $\sgn_{\rmO^\pm_{2n}}\cdot R_x$ and $\chi_{\rmO^\pm_{2n}} R_x$
are almost characters of\/ $\rmO^\pm_{2n}$.
\end{lem}
\begin{proof}
The proof is the same as that of Lemma~\ref{0422}.
\end{proof}

Then we define the $\Theta$-rank of an almost character $R_x\in\cala(\bfG)$
for $\bfG=\Sp_{2n}$, $\rmO^\pm_{2n}$, or $\SO_{2n+1}$ as follows.
\begin{itemize}
\item For $R_x\in\cala(\Sp_{2n})$,
we define
\[
\Theta\text{\rm -rk}(R_x)
=\min\{\,2k,2k+1\mid R_x\text{ occurs in }\widetilde\Theta^\psi_{\Sp_{2n},\rmO^\pm_{2k}}
\text{ or }\widetilde\Theta^\psi_{\Sp_{2n},\SO_{2k+1}}\text{ for some }\psi\,\};
\]

\item for $R_x\in\cala(\rmO^\pm_{2n})$,
we define
\[
\Theta\text{\rm -rk}(R_x)
=\min\{\,2k\mid\text{one of }R_x,R_{\bar x},R_{x^\rmt},R_{\bar x^\rmt}
\text{ occurs in }\widetilde\Theta^\psi_{\rmO^\pm_{2n},\Sp_{2k}}\text{ for some }\psi\,\},
\]
where $\bar x^\rmt=(\bar x)^\rmt$;

\item for $R_x\in\cala(\SO_{2n+1})$,
we define
\[
\Theta\text{\rm -rk}(R_x)
=\min\{\,2k\mid\text{one of }R_x,R_{\bar x}
\text{ occurs in }\widetilde\Theta^\psi_{\SO_{2n+1},\Sp_{2k}}\text{ for some }\psi\,\}.
\]
\end{itemize}
Note that $\bar x$ is also defined for $x\in X(\SO_{2n+1})$ in Lemma~\ref{0422}.

\begin{cor}
Let $\bfG=\Sp_{2n}$, $\rmO^\pm_{2n}$, or $\SO_{2n+1}$.
For any $x\in X(\bfG)$,
we have
\[
\Theta\text{\rm -rk}(\rho_x)=\Theta\text{\rm -rk}(R_x).
\]
\end{cor}
\begin{proof}
This follows immediately from the definition of $\Theta\text{\rm -rk}(R_x)$ above,
the definition of $\Theta\text{\rm -rk}(\rho_x)$ in \cite{pan-theta-rank}, Theorem~\ref{0101}.
\end{proof} 

\bibliography{refer}
\bibliographystyle{amsalpha}

\end{document}